\documentclass[11pt]{amsart}
\usepackage{amssymb}
\usepackage{amsmath,amssymb}

\theoremstyle{plain}
\newtheorem{thm}{Theorem}[section]
\newtheorem{theorem}[thm]{Theorem}

\newtheorem{lemma}[thm]{Lemma}
\newtheorem{corollary}[thm]{Corollary}
\newtheorem{proposition}[thm]{Proposition}
\theoremstyle{definition}
\newtheorem{remark}[thm]{Remark}

\newtheorem{notation}[thm]{Notation}

\newtheorem{definition}[thm]{Definition}

\newtheorem{assumption}[thm]{Assumption}

\newtheorem{example}[thm]{Example}

\newtheorem{setup}[thm]{Setup}
\numberwithin{equation}{section}
\newcommand{\II}{{\rm II}}
\newcommand{\III}{{\rm III}}

\newcommand{\sA}{{\mathcal A}}

\newcommand{\sC}{{\mathcal C}}
\newcommand{\sD}{{\mathcal D}}

\newcommand{\sG}{{\mathcal G}}

\newcommand{\sK}{{\mathcal K}}
\newcommand{\sL}{{\mathcal L}}

\newcommand{\sN}{{\mathcal N}}
\newcommand{\sO}{{\mathcal O}}
\newcommand{\sP}{{\mathcal P}}
\newcommand{\sQ}{{\mathcal Q}}
\newcommand{\sR}{{\mathcal R}}
\newcommand{\sS}{{\mathcal S}}
\newcommand{\sT}{{\mathcal T}}
\newcommand{\sU}{{\mathcal U}}
\newcommand{\sV}{{\mathcal V}}
\newcommand{\sW}{{\mathcal W}}
\newcommand{\sX}{{\mathcal X}}
\newcommand{\sY}{{\mathcal Y}}
\newcommand{\sZ}{{\mathcal Z}}

\newcommand{\C}{{\mathbb C}}

\newcommand{\N}{{\mathbb N}}
\newcommand{\BP}{{\mathbb P}}

\newcommand{\Q}{{\mathbb Q}}

\newcommand{\BS}{{\mathbb S}}

\newcommand{\Z}{{\mathbb Z}}

\newcommand{\rH}{{\rm Hom}}

\newcommand{\End}{{\rm End}}

\newcommand{\bL}{\textbf{L}}
\newcommand{\bG}{\textbf{G}}
\newcommand{\bX}{\textbf{X}}
\newcommand{\bl}{\textbf{l}}

\newcommand{\bv}{\textbf{v}}

\newcommand{\fg}{{\mathfrak g}}
\newcommand{\fsl}{{\mathfrak s}{\mathfrak l}}
\newcommand{\fgl}{{\mathfrak g}{\mathfrak l}}

\newcommand{\fl}{{\mathfrak l}}

\newcommand{\aut}{{\mathfrak a}{\mathfrak u}{\mathfrak t}}

\newcommand\Aut{\rm Aut}

\def\Sym{\mathop{\rm Sym}\nolimits}
\def\Pic{\mathop{\rm Pic}\nolimits}
\def\Hom{\mathop{\rm Hom}\nolimits}

\title[${\rm G}_2$-horospherical manifold]{Recognizing the ${\rm G}_2$-horospherical manifold of Picard number 1 by varieties of minimal rational tangents}

\author{Jun-Muk Hwang and Qifeng Li}

\thanks{This work was supported by the Institute for Basic Science (IBS-R032-D1).}

\begin{document}

\begin{abstract}
Pasquier and Perrin discovered that the ${\rm G}_2$-horospherical manifold ${\bf X}$ of Picard number 1  can be realized as a smooth specialization of the rational homogeneous space parameterizing the lines on the 5-dimensional hyperquadric, in other words, it can be deformed nontrivially to the rational homogeneous space. We show that ${\bf X}$ is the only smooth projective variety with this property. This is obtained as a consequence of our main result that ${\bf X}$ can be recognized by its VMRT, namely,
a Fano manifold of Picard number 1 is biregular to ${\bf X}$ if and only if its VMRT at a general point is projectively isomorphic to that of ${\bf X}$.
We employ the method the authors developed to solve the corresponding problem  for symplectic Grassmannians, which constructs a flat Cartan connection in a neighborhood of a general minimal rational curve. In adapting this method to ${\bf X}$, we need an intricate study of the positivity/negativity of vector bundles with respect to a family of rational curves, which is subtler than the case of symplectic Grassmannians because of the nature of the differential geometric structure on ${\bf X}$ arising from VMRT.
\end{abstract}

\maketitle


\section{Introduction}\label{s.I}

We work in the holomorphic category: all varieties, groups and algebras are defined over complex numbers and all morphisms are holomorphic.  Unless mentioned otherwise, open sets refer to
Euclidean topology.  The projectivization $\BP
V$ of a vector space (or a vector bundle) $V$ is taken in the
classical sense, i.e., as the set of 1-dimensional subspaces of $V$.

We are mainly interested in the ${\rm G}_2$-horospherical manifold of Picard number 1, a 7-dimensional Fano manifold $\bX$ defined as follows. Let $\bL$ be a  simple Lie group of type ${\rm G}_2$ with the short fundamental weight $\omega_1$ and the long fundamental weight $\omega_2$.   Let $V(\omega_1)$ (resp. $V(\omega_2)$) be an irreducible representation of $\bL$ with the highest weight $\omega_1$ (resp. $\omega_2$) and let $v_1 \in V(\omega_1)$ (resp. $v_2 \in V(\omega_2)$ ) be a highest weight vector.
Define ${\bf X} \subset  \BP(V(\omega_1) \oplus V(\omega_2))$ as the closure of the orbit of $\bL$ through $[v_1 + v_2] \in \BP(V(\omega_1) \oplus V(\omega_2)).$
It was proved by Pasquier that $\bX$ is a smooth projective variety of Picard number 1. In fact, it is one of the smooth horospherical varieties of Picard number 1 classified in Theorem 0.1 of \cite{Pa}.  Moreover, the identity component $\bG$ of the automorphism group of $\bX$ is isogenous to $(\bL \times \C^*) \ltimes  V(\omega_1)$
and it has two orbits in $\bX$ (Theorem 1.11 of \cite{Pa}).

In \cite{PP}, Pasquier and Perrin discovered a  remarkable deformation property of $\bX.$
To explain it, let $\Q^n \subset \BP^{n+1}, n \geq 5,$ be the $n$-dimensional smooth quadric hypersurface and let ${\rm Lines}(\Q^n)$ be the space of lines lying on $\Q^n$.
Both $\Q^n$ and ${\rm Lines}(\Q^n)$ are  smooth projective varieties of
Picard number 1, homogeneous under the group ${\bf
PSO}(n+2)$ of  automorphisms of  $\Q^n$.
In \cite{Hw97} and also in \cite{HM02},  the
following was proved for $n \geq 6$.

\begin{theorem}\label{t.97} Let $\pi: \sX \to \Delta$ be a smooth projective
morphism from a complex manifold $\sX$ to the unit disc $\Delta
\subset \C$. If the fiber $\pi^{-1}(t)$  is biregular to ${\rm Lines}(\Q^n)$
for each $t\in \Delta \setminus \{0\}$, then the central  fiber
$\pi^{-1}(0)$ is biregular to ${\rm Lines}(\Q^n)$, too.
\end{theorem}

In the statement of Theorem \ref{t.97} given in \cite{Hw97} and \cite{HM02},  the case of $n = 5$  was included erroneously, although neither the proof in \cite{Hw97} nor the slightly
different proof in \cite{HM02} work when $n=5$. The discovery of
Pasquier and Perrin (Proposition 2.3 of \cite{PP}) is the following   counter-example to the
statement of Theorem \ref{t.97} when $n=5$.

\begin{theorem}\label{t.PP}
There exists a smooth projective
morphism  $\pi: \sX \to \Delta$ from an 8-dimensional complex manifold $\sX$ to the unit disc $\Delta
\subset \C$ such that the central  fiber
$\pi^{-1}(0)$ is biregular to $\bX$ while the  fiber $\pi^{-1}(t)$  for each $t\in \Delta \setminus \{0\}$ is biregular to ${\rm Lines}(\Q^5)$. \end{theorem}

Why does the proof of Theorem \ref{t.97} in \cite{Hw97} and \cite{HM02} fail when $n=5$?
It turns out that the key issue  lies in the description of
the {\em embedded deformations} of the Segre embedding $\BP^1 \times
\Q^{n-4} \subset \BP^{2n-5}.$ To be precise,  consider a  family
of smooth projective varieties $$\{ \sC_t \subset \BP^{2n-5},  t
\in \Delta \}$$ such that for $t \neq 0$, the submanifold $\sC_t
\subset \BP^{2n-5}$ is isomorphic to the Segre embedding $\BP^1 \times
\Q^{n-4} \subset \BP^{2n-5}.$ The question is
whether the embedding $\sC_0 \subset \BP^{2n-5}$ is  isomorphic to
the Segre embedding, too. It is not hard to see that this is so if $n
\geq 6$ (two slightly different proofs are given in Lemma 4 of
\cite{Hw97} and Lemma 3 of \cite{HM02}).  But when $n=5$, one additional possibility of
$\sC_0 \subset \BP^{2n-5}$ exists. Let $\BS$ be the Hirzebruch
surface $\BP(\sO(-1) \oplus \sO(-3))$ where $\sO(k)$ denotes the
line bundle on $\BP^1$ of degree $k$. It can be
embedded in the projective space $$\BP H^0(\BP^1, \sO(3) \oplus
\sO(1))^* \cong \BP^5.$$ There exists a family of
smooth surfaces $$\{ \sC_t \subset \BP^5, t \in \Delta\}$$ such that the submanifold $\sC_t \subset \BP^5$ is the Segre embedding of
$\BP^1 \times \Q^1$ for $t \neq 0$, but $\sC_0 \subset \BP^5$ is
isomorphic to $\BS \subset \BP^5$. This additional possibility was
overlooked in \cite{Hw97} and \cite{HM02}.

In the setting of Theorem \ref{t.97}, the above situation arises in
the following way (see the proof of Lemma \ref{l.deform} for more details).
Take a section $\sigma: \Delta \to \sX$ of $\pi$ such that
$\sigma(0)$ is a general point of the central fiber. Then consider
the family $ \nu: \sK_{\sigma} \to \Delta$ whose fiber $\sK_{\sigma(t)}$ at $t \in \Delta$
is the (normalized) Chow space of minimal rational
curves on $\pi^{-1}(t)$ passing through $\sigma(t)$. It is known
that $\nu$ is a smooth projective morphism and the general fiber
$\nu^{-1}(t)$ is biregular to $\BP^1 \times \Q^{n-4}.$ In the proofs of Theorem \ref{t.97}
in \cite{Hw97} and \cite{HM02}, the first step is to show that
$\nu^{-1}(0)$ is also biregular to $\BP^1 \times \Q^{n-4}$, which implies that the varieties of minimal rational tangents (to be abbreviated VMRT)  at general points of $\pi^{-1}(0)$ are isomorphic to the Segre embedding of $\BP^1 \times \Q^{n-4}$ (see  Definition \ref{d.VMRT} for the definition of VMRT). This step is essentially equal to the embedded
deformation problem for surfaces discussed in the previous paragraph.

If this step works, then the rest of the arguments  in \cite{Hw97} and \cite{HM02} works perfectly well even when $n=5$.  The rest of the arguments, deducing $\pi^{-1}(0) \cong {\rm Lines}(\Q^n)$ from the isomorphism type of the VMRT, has been developed into  the following stronger result by N. Mok in \cite{Mk}.

\begin{theorem}\label{t.Mok}
Let $X$ be a Fano manifold of Picard number 1 whose
VMRT at a general point is
isomorphic (as projective subvarieties of the projectivized tangent space)
to that of ${\rm Lines}(\Q^n),  n \geq 5.$ Then $X$ is biregular
to ${\rm Lines}(\Q^n)$.
\end{theorem}

In \cite{Mk}, Mok was suggesting that problems like Theorem \ref{t.Mok},
asking to recognize a Fano manifold from its VMRT, are more fundamental than
deformation rigidity like Theorem \ref{t.97}. The discovery of Theorem \ref{t.PP} supports this viewpoint. It is natural to expect that a version of Theorem \ref{t.Mok} should hold when ${\rm Lines}(\Q^n)$ is replaced by
any rational homogeneous space $G/P$ of Picard number 1. Indeed, this expectation was verified for any $G/P$ associated to a long simple root in \cite{Mk} and \cite{HH}.
The case of $G/P$ associated to a short simple root turns out to be much harder and
only recently, the authors have been able to verify it for symplectic Grassmannians in \cite{HL} by developing a novel method of constructing Cartan connections.

The goal of this paper is to use the method developed in \cite{HL} to prove the following version of Theorem \ref{t.Mok} for $\bX$.

\begin{theorem}\label{t.main}
Let $X$ be a 7-dimensional Fano manifold of Picard number 1 whose
VMRT at a general point is
isomorphic (as projective subvarieties of the projectivized tangent space)
to that of $\bX.$ Then $X$ is biregular
to $\bX$.
\end{theorem}

As applications of Theorem \ref{t.main}, we prove the following two results on deformation.

\begin{theorem}\label{t.deform}
Let $\pi: \sX \to \Delta$ be a smooth projective morphism from an
8-dimensional complex manifold $\sX$ to the unit disc $\Delta
\subset \C$. If the fiber $\pi^{-1}(t)$  is biregular to ${\rm Lines}(\Q^5)$
for each $t\in \Delta \setminus \{0\}$, then the central fiber
$\pi^{-1}(0)$ is biregular to either ${\rm Lines}(\Q^5)$ or  $\bX$. \end{theorem}

\begin{theorem}\label{t.rigid}
Let $\pi: \sX \to \Delta$ be a smooth projective morphism from an
8-dimensional complex manifold $\sX$ to the unit disc $\Delta
\subset \C$. If the fiber $\pi^{-1}(t)$  is biregular to  $\bX$  for each $t\in
\Delta \setminus \{0\}$, then the central fiber $\pi^{-1}(0)$ is
biregular to $\bX$, too.  \end{theorem}

Theorem \ref{t.deform} is the correct version of Theorem \ref{t.97} when $n=5$.
Notice the difference between Theorem \ref{t.97} (or Theorem \ref{t.rigid}) and Theorem \ref{t.deform}.
While Theorem \ref{t.97} is a rigidity result stating the {\em continuity} of certain structures, Theorem
\ref{t.deform}  is to {\em determine} a new structure popping up at the limit. Although Theorem \ref{t.97} could be proved without using Theorem \ref{t.Mok} as the proofs in \cite{Hw97} and \cite{HM02} show,
  Theorem \ref{t.main} is essential to prove Theorem \ref{t.deform}.

The proof of Theorem \ref{t.main} follows the strategy of \cite{HL}. In particular, it uses the differential geometric machinery developed in Section 2 of \cite{HL}, which
we recall  in Section \ref{s.HL} of this paper. The key result   is Theorem \ref{t.HL} below which gives   a way to  recognize  homogenous spaces by checking two conditions:
({\sf i}) the prolongation property of the associated graded Lie algebra and ({\sf ii}) vanishing of sections of certain cohomological vector bundles.   In Section \ref{s.prolongation}, we verify the prolongation property ({\sf i}) for the Lie algebra $\fg$ of the automorphism group of $\bX$. In Section \ref{s.S}, we review some projective-geometric properties of the surface $\BS$, which are needed to study the geometric structures modeled on $\bX$.
In Section \ref{s.vanishing}, we check the vanishing condition ({\sf ii}) for the geometric structure modeled on $\bX$.   After recalling  some basic results on VMRT in Section \ref{s.vmrt},  we prove Theorems \ref{t.main}, Theorem  \ref{t.deform}  and Theorem \ref{t.rigid} in Section \ref{s.proof}. The overall line of arguments is parallel to that of \cite{HL}, although the details are different because of the difference in the nature of the underlying differential geometric objects. Especially, to check the vanishing condition in Section \ref{s.vanishing}, some new notions of positivity/negativity of vector bundles along a given family of rational curves are introduced in Definition \ref{d.generation}, which are  more intricate than the SAF-condition used in Section 5 of \cite{HL}.

It seems worthy to add a remark on the history of this paper.  Several years ago, the first-named author had a preprint titled `Deformation of the space of lines on the 5-dimensional hyperquadric' where an approach to Theorems \ref{t.main} was presented in a long and complicated  argument.
The preprint had been circulated among some colleagues, but had never been published because some details of the argument needed clarification. The current paper has borrowed some part from that preprint,  but its key methodology based on  \cite{HL} is completely different and conceptually simpler.

\section{Review of Section 2 of \cite{HL}}\label{s.HL}
In this section, we briefly recall the results of Section 2 in \cite{HL}.

\begin{definition}\label{d.filtration}
Let $M$ be a complex manifold. A {\em filtration } $F^{\bullet}$ on $M$ is a collection of subbundles
$(F^j, j \in \Z)$ of the tangent bundle $TM$ such that
 \begin{itemize}
 \item[(i)]
 $F^{j+1} \subset F^{j}$ is a vector subbundle for each $j \in \Z$;
 \item[(ii)] $F^{-\nu} = TM$ for some $\nu \in \N$;
 and
 \item[(iii)] regarding the locally free sheaf $\sO(F^i)$ as a sheaf of vector fields on $M$, we have $$[\sO(F^{-i}), \sO( F^{-j})] \subset \sO(F^{-k})$$ for the Lie bracket of vector fields and for any nonnegative integers
     $i, j, k,$
 satisfying $i+j \leq k$. \end{itemize}
 For a point $x \in M$, we will denote by $F^{k}_x $ the fiber of $F^k $ at the point $x$. By (iii), the graded vector space $${\rm Symb}_x(F^{\bullet}) \ := \ \bigoplus_{i \in \N} F_x^{-i}/F_x^{-i+1} $$
 for each point $x \in M$ has the structure of a nilpotent graded Lie algebra, called the {\em symbol algebra} of the filtration.
 \end{definition}

 \begin{definition}\label{d.G_0-structure}
For a graded nilpotent Lie algebra $\fg_{-}= \fg_{-1} \oplus \cdots \oplus \fg_{-\nu},$
a {\em filtration of type} $\fg_{-}$ on a complex manifold $M$ is a filtration $F^{\bullet} = (F^j, j \in \Z)$ on $M$ such that
\begin{itemize}
\item[(i)] $F^k =0$ for all $k \geq 0$; \item[(ii)] $F^{-k} M = TM$ for all $k \geq \nu$; and
\item[(iii)] for any $x\in M$, the symbol algebra $${\rm Symb}_x(F^{\bullet}) = \bigoplus_{i \in \N} F^{-i}_{x}/F^{-i+1}_{x}$$  is isomorphic to $\fg_{-}$ as graded Lie algebras. \end{itemize} Let ${\rm grAut}(\fg_-)$ be the group of Lie algebra automorphisms of $\fg_-$ preserving the gradation.
The {\em graded frame bundle} of the manifold $M$ with a filtration $F^{\bullet}$ of type $\fg_{-}$  is the  ${\rm gr}\Aut(\fg_{-})$-principal bundle ${\rm grFr}(F^{\bullet})$ on $M$ whose fiber at $x$ is the set of graded Lie algebra isomorphisms from
$\fg_{-}$ to ${\rm Symb}_x(F^{\bullet})$. Let $G_0 \subset {\rm gr}\Aut(\fg_{-})$ be a connected algebraic subgroup.  A $G_0$-{\em structure subordinate to the filtration} $F^{\bullet}$ on $M$ means a $G_0$-principal subbundle  $\sA \subset {\rm grFr}(F^{\bullet})$.\end{definition}

\begin{definition}\label{d.prolongation}
Let $\fg_-= \fg_{-1} \oplus \cdots \oplus \fg_{-\nu}$ be a graded nilpotent Lie algebra. Let ${\rm gr}\aut(\fg_-)$ be the  Lie algebra of ${\rm grAut}(\fg_-)$. Fix a connected algebraic subgroup $G_0 \subset {\rm grAut}(\fg_-)$ and its Lie algebra $\fg_0 \subset {\rm gr} \aut(\fg_-)$. The \emph{universal prolongation} of $(\fg_0, \fg_-)$ is a graded Lie algebra $\oplus_{k \geq - \nu} \fg_k$ whose adjoint action on $\fg_-$ is effective and is maximal among all graded Lie algebra extending $\fg_0 \oplus \fg_-$ whose adjoint action on $\fg_-$ is effective. The universal prolongation exists for any graded nilpotent Lie algebra and is unique up to isomorphisms. See Definition 2.7  of \cite{HL} or Section 2.3 of \cite{Ya} (where it is called `algebraic prolongation') for an explicit description of the universal prolongation. \end{definition}

\begin{definition}\label{d.cohomology}
Let $\oplus_{k \geq - \nu} \fg_k$ be the universal prolongation of $(\fg_0, \fg_-)$ in Definition \ref{d.prolongation}.
Define, for each  $\ell \in \N$,
\begin{eqnarray*}
C^{\ell, 1}(\mathfrak{g}_{0}) & := & \bigoplus_{j\in \N}\Hom(\fg_{-j}, \mathfrak{g}_{-j+\ell}), \\
C^{\ell, 2}(\fg_0) & :=&
\bigoplus_{i, j \in \N} \Hom(\fg_{-i} \wedge \fg_{-j}, \fg_{- i -j + \ell}).
\end{eqnarray*}
For an element $f \in C^{\ell, 1}(\fg_0)$, define $\partial f \in C^{\ell, 2}(\mathfrak{g}_0)$ by
$$\partial f(u, v)  = [f(u), v] + [u, f(v)] -f([u, v]) \mbox{ for all } u, v \in \fg_{-}.$$ This determines a $G_0$-module homomorphism $$\partial: C^{\ell, 1}(\fg_0) \to C^{\ell, 2}(\fg_0).$$  Its cokernel $C^{\ell, 2}(\fg_0)/\partial(C^{\ell,1}(\fg_0))$ has a natural $G_0$-module structure. \end{definition}

\begin{example}\label{e.model}
In Definition \ref{d.cohomology}, assume that the universal prolongation $\fg:= \oplus_{k \geq - \nu} \fg_k$ is finite-dimensional.
Let $G$ be a connected Lie group with Lie algebra $\fg$ and $G^0 \subset G$ be the closed subgroup with Lie algebra $\oplus_{k \geq 0} \fg_k$.
Then the homogenous space $G/G^0$ has a natural filtration $F^{\bullet}_{G/G^0}$ of type $\fg_-$ and a $G_0$-structure $\sA_{G/G^0} \subset {\rm grFr}(F^{\bullet}_{G/G^0}) $ subordinate to the filtration. \end{example}

The following is the main result in Section 2 of \cite{HL}: it is Theorem 2.17 combined with Proposition 2.5 in \cite{HL}.

\begin{theorem}\label{t.HL}
In Definition \ref{d.cohomology}, assume that the action of $G_0$ on $\fg_-$ has no nonzero fixed vector and the universal prolongation $\fg:= \oplus_{k \geq - \nu} \fg_k$ is finite-dimensional. Let $G$ be a connected Lie group of adjoint type with Lie algebra $\fg$ and $G^0 \subset G$ be the connected closed subgroup with Lie algebra $\oplus_{k \geq 0} \fg_k$. Let $M$ be a simply connected complex manifold  with a filtration $F^{\bullet}$ of type $\fg_{-}$ and let $\sA \subset {\rm grFr}(F^{\bullet})$ be a $G_0$-structure on $M$ subordinate to the filtration. Let $$C^{\ell,2}(\sA)/\partial(C^{\ell,1}(\sA))$$ be the vector bundle associated to the principal $G_0$-bundle $\sA$ by the natural representation of $G_0$ on
$$C^{\ell, 2}(\fg_0)/\partial(C^{\ell,1}(\fg_0))$$ from Definition \ref{d.cohomology}.  Assume that $H^0(M, C^{\ell,2}(\sA)/\partial(C^{\ell,1}(\sA)) ) =0$ for all $ \ell \geq 1.$ Then there exists an open immersion $\phi: M \to G/G^0$ that sends $F^{\bullet}$ and  $\sA$ isomorphically to $ F^{\bullet}_{G/G^0}|_{\phi(M)}$ and $\sA_{G/G^0}|_{\phi(M)}$ in Example \ref{e.model}.
\end{theorem}

The statement \ref{t.HL} is slightly different from Theorem 2.17 of \cite{HL}. The latter theorem gives the existence of a Cartan connection on a prolonged principal bundle under the weaker assumption of the vanishing $$H^0(M, C^{\ell,2}(\sA)/\partial(C^{\ell,1}(\sA)) ) =0 \mbox{ for all } 1 \leq \ell \leq\mu + \nu.$$ By assuming the vanishing for all $\ell \geq 1$,  we can obtain the local flatness of the Cartan connection by the same argument as in the proof of Theorem 2.17 of \cite{HL}, because the curvature of the resulting Cartan connection would give elements of $H^0(M, C^{\ell,2}(\sA)/\partial(C^{\ell,1}(\sA)) )$ for some values of $\ell \geq 1$. This implies the existence of the open immersion $f$ by Proposition 2.5 of \cite{HL}.

\section{Structure of the graded Lie algebra $\fg$}\label{s.prolongation}

In this section, we introduce a graded Lie algebra structure on the Lie algebra $\fg$ of the automorphism group of the variety $\bX$ in Section \ref{s.I}  and examine its prolongation property.

To start with, we recall the notation and some  convention on tensor products of vector spaces.
Let $U$ be a vector space of dimension $n$, and $U^*$ be its dual vector space. For each $h\in U^*$ and $u\in U$, we regard $h\otimes u$ as the unique element in $\fgl(U)$ sending each $v\in U$ to $h(v)u\in U$. This gives an identification of vector spaces $U^*\otimes U=\fgl(U)$, and we can transfer the Lie algebra structure of $\fgl(U)$ to $U^*\otimes U$ via this identification.

We follow the  convention of  Appendix B of \cite{FH}  for the exterior and symmetric powers. In particular, for nonnegative integers $p$ and $q$,  we have natural pairings
\begin{eqnarray*}
\langle\, , \, \rangle: \wedge^pU\otimes\wedge^pU^*\rightarrow\C, & \langle\, , \, \rangle: \Sym^pU\otimes\Sym^pU^*\rightarrow\C,
\end{eqnarray*}
and natural contractions
\begin{eqnarray*}
\lrcorner: \wedge^p U\otimes\wedge^{p+q} U^*\rightarrow\wedge^q U^*, & \lrcorner: \Sym^p U\otimes\Sym^{p+q} U^*\rightarrow\Sym^q U^*,\\
\llcorner: \wedge^{p+q}U\otimes\wedge^p U^*\rightarrow\wedge^q U, & \llcorner: \Sym^{p+q}U\otimes\Sym^p U^*\rightarrow\Sym^q U,
\end{eqnarray*}
such that the following holds for
 any $h, h_1, h_2\in U^*$ and any $u, u_1, u_2\in U$.
\begin{eqnarray*}
 \langle u, h\rangle &=& u\llcorner h=u\lrcorner h=h(u), \\
 (u_1\wedge u_2)\llcorner h&=&h(u_1)u_2-h(u_2)u_1, \\
 u\lrcorner(h_1\wedge h_2)&=&h_2(u)h_1-h_1(u)h_2, \\
 \langle u_1\wedge u_2, h_1\wedge h_2\rangle &=& h_1(u_1)h_2(u_2)-h_1(u_2)h_2(u_1), \\
 u^q\llcorner h&=&qh(u) u^{q-1}, \\
 u\lrcorner h^q&=&qh(u) h^{q-1}, \\
 u^q\llcorner h^p&=&(u^q\llcorner h)\llcorner h^{p-1}, \\
 u^p\lrcorner h^q&=&u^{p-1}\lrcorner(u\lrcorner h^q), \\
 \langle u^q, h^q\rangle = u^q\llcorner h^q &=& u^q\lrcorner h^q=q!h(u)^q.
\end{eqnarray*}
One can also see that
$$\langle u_1\wedge u_2, h_1\wedge h_2\rangle=((u_1\wedge u_2)\llcorner h_1)\llcorner h_2=u_1\lrcorner(u_2\lrcorner(h_1\wedge h_2)).\\$$

Now, let $\fl$ be the simple Lie algebra of type ${\rm G}_2$ with a system of simple roots $\{ \alpha_1, \alpha_2\}$ with respect to a Cartan subalgebra. We use the convention that $\alpha_1$ is short and $\alpha_2$ is long such that their coroots satisfy
$$\alpha_1 (\alpha_2^{\vee}) = -1, \ \  \alpha_2(\alpha_1^{\vee}) = -3.$$ Then the fundamental weights are
$$ \omega_1 = 2 \alpha_1 + \alpha_2, \ \ \omega_2 = 3 \alpha_1 + 2 \alpha_2.$$
We consider the gradation $$\fl = \fl_{-2} \oplus \fl_{-1}  \oplus \fl_0 \oplus \fl_1  \oplus \fl_2$$ determined by the coefficients of the long root $\alpha_2$. The subalgebra $\fl^0 := \fl_0 \oplus \fl_1  \oplus \fl_2$ is a parabolic subalgebra of $\fl$ and there is a unique element $E \in \fl_0$, called the characteristic element of the gradation, which acts on $\fl_i$ with eigenvalue  $i$ under the adjoint representation (Section 3.1 of \cite{Ya}). It is easy to check (e.g. from  the  partition of the set of positive roots in p. 446 of \cite{Ya}) that $$E= \alpha_1^{\vee} + 2 \alpha_2^{\vee}, \ \ \omega_1(E) = 1, \ \ \omega_2(E) = 2.$$
Then $$\fl_0 = \C E \oplus  [\fl_0, \fl_0] = \C E \oplus (\fl_{-\alpha_1} + \C \alpha_1^{\vee} + \fl_{\alpha_1}) $$ and the Lie subalgebra $[\fl_0, \fl_0]$ is isomorphic to $\fsl_2$.

Let $V$ be the irreducible representation of $\fl$ with the highest weight $\omega_1$. The 7-dimensional vector space $V$ has weights (e.g. the second diagram in Section 22.3 of \cite{FH})
$$0, \pm \omega_1, \pm( \omega_1-\omega_2), \pm (2 \omega_1 - \omega_2).$$
We have the $E$-eigenspace decomposition $V = V_{-1}  \oplus V_0  \oplus V_1$ with weights
$$\{ - \omega_1, \omega_1- \omega_2\}, \ \{ 0, \pm(2 \omega_1 - \omega_2)\}, \ \{ \omega_1, \omega_2-\omega_1\}$$ for $V_{-1}, V_0, V_1, $ respectively.
It is easy to check that each $V_i, i =-1, 0, 1,$ is an irreducible representation of $[\fl_0, \fl_0] \cong \fsl_2$.

\begin{definition}\label{d.fg}
Let $\fg = \fg_{-2}  \oplus \fg_{-1}  \oplus \fg_0  \oplus \fg_1  \oplus \fg_2$ be the semi-direct product $(\fl \oplus \C{\rm Id}_V) \ltimes V$ with the gradation \begin{eqnarray*}
\fg_{-2} & :=& \fl_{-2} \\
\fg_{-1} & := & \fl_{-1} \oplus V_{-1} \\
\fg_0 & := & \fl_0 \oplus \C{\rm Id}_V \oplus V_0 \\
\fg_1 & := & \fl_1 \oplus V_1 \\
\fg_2 & :=& \fl_2. \end{eqnarray*} \end{definition}

The following lemma is immediate.

\begin{lemma}\label{l.filt}
For each $k \in \Z$, the filtration $0 \subset V_k \subset \fg_k$ is a $\fg_0$-reductive filtration in the sense of Definition 5.2 of \cite{HL},  i.e., the subspace $V_k$ is a $\fg_0$-submodule of $\fg_k$, and the nil-radical of $\fg_0$ sends $V_k$ to 0 and $\fg_k$ to $V_k$ for each $k$ such that the graded objects of the filtration become representations of the quotient Lie algebra of $\fg_0$ by its nil-radical.  \end{lemma}

Recall the following result from Proposition 49 of \cite{Kim}.

\begin{proposition}\label{p.kim}
The graded Lie algebra $\fg$  in Definition \ref{d.fg} is the universal prolongation of $(\fg_0, \fg_-:= \fg_{-1} \oplus \fg_{-2} )$ in the sense of Definition \ref{d.prolongation}. \end{proposition}

We can describe the graded Lie algebra structure of  $\fg$ as follows.

\begin{proposition}\label{p.repre}
Fix  vector spaces $W$ and  $R$ with $\dim W =2$ and $\dim R =1$.  Putting $Q:= R^{\otimes 2} \otimes (\wedge^2 W)^{\otimes 3}$,  define
\begin{eqnarray*}\bl_{-2} & = & Q \\
\bl_{-1} & = & R \otimes \Sym^3 W \\
\bl_1 & = & R^* \otimes \Sym^3 W^* \\
\bl_2 & = & Q^* \\
\bv_{-1} &=&  W \\
\bv_0 &=& R^* \otimes \Sym^2 W^* \\
\bv_1 & =& Q^* \otimes W. \end{eqnarray*}
Define a  graded nilpotent Lie algebra structure on $\bl_{-1} \oplus \bl_{-2}$  by
$$ [r \otimes w_1^3, r \otimes w_2^3] = r^2 \otimes (w_1 \wedge w_2)^3 $$ for any $r \in R$ and $w_1, w_2 \in W.$
Then the following holds.
\begin{itemize} \item[(1)]
The natural representation of the Lie algebra $\fgl(R) \oplus \fgl(W)$  on
$$\bl_{-1} \oplus \bl_{-2} = (R \otimes \Sym^3 W)  \oplus Q$$ respects the structure of the graded Lie algebra.
The kernel of this representation is $$ \C({\rm Id}_W - 3 {\rm Id}_R) \subset \fgl(R) \oplus \fgl(W).$$
Under the natural inclusion of Lie algebras $\fgl(R)\oplus\fgl(W)\subset\fgl(R\oplus W)$,  consider the Lie algebra
$$\bl_0=(\fgl(R)\oplus\fgl(W))\cap\fsl(R\oplus W).$$
Then  $\C({\rm Id}_W-2{\rm Id}_R)$ is the center of $\bl_0$, and
\begin{eqnarray*}
 \bl_0 &=&\fsl(W)\oplus\C({\rm Id}_W-2{\rm Id}_R),\\
\fgl(R)\oplus\fgl(W) &=& \bl_0\oplus\C({\rm Id}_W - 3 {\rm Id}_R).
\end{eqnarray*}
The natural representation of $\fgl(R) \oplus \fgl(W)$ on each of the vector spaces $\bl_1, \bl_2, \bv_{-1}, \bv_0$ and $\bv_1$ induces a representation of $\bl_0$ on it.
\item[(2)] Fix a nonzero $r \in R$ and its dual $r^* \in R^*$.
Consider  the following natural operations
for any $w_1,  w_2 \in W $ and $h_1, h_2 \in W^*$.
\begin{itemize} \item[$\bullet$] $[\bl_1, \bl_{-2}]$:  $$ [ r^*\otimes h_1^3, r^2 \otimes (w_1 \wedge w_2)^3] = -r
\otimes((w_1\wedge w_2)\llcorner h_1)^3\in\bl_{-1}.$$
\item[$\bullet$]
$[\bl_1, \bl_{-1}]$: $$ [r^* \otimes h_1^3, r \otimes w_1^3] = -h_1(w_1)^3 {\rm Id}_R + h_1(w_1)^2 h_1\otimes w_1\in\bl_0. $$
\item[$\bullet$]
$[\bl_1, \bl_1]$:  $$ [ r^* \otimes h_1^3, r^* \otimes h_2^3] = (r^*)^2 \otimes (h_1 \wedge h_2)^3\in\bl_2,$$
\item[$\bullet$]
$[\bl_2, \bl_{-2}]$:
\begin{align*}
& [ (r^*)^2 \otimes (h_1 \wedge h_2)^3, r^2 \otimes (w_1 \wedge w_2)^3 ] \\
=& \langle w_1\wedge w_2, h_1\wedge h_2\rangle^3 (-2 {\rm Id}_R + {\rm Id}_W)\in\bl_0,
\end{align*}
\item[$\bullet$]
$[\bl_2, \bl_{-1}]$: $$ [ (r^*)^2 \otimes (h_1 \wedge h_2)^3, r \otimes w_1^3] = - r^* \otimes(w_1\lrcorner(h_1\wedge h_2))^3\in\bl_1.$$
\end{itemize}
Together with the natural operation of $\bl_0$ on $\bl_i, -2 \leq i \leq 2$,
the above operations define a graded Lie algebra structure on $\bl := \oplus_{i=-2}^2 \bl_i.$
\item[(3)] Fix a nonzero $ f \in \wedge^2 W$ and its dual $f^* \in \wedge^2 W^*$. Consider the following operations for any $w_1,  w_2, w_3 \in W $ and $h_1, h_2, h_3 \in W^*$,
    \begin{itemize}
\item[$\bullet$]
$[\bl_2, \bv_{-1}]$:
\begin{align*}
& [(r^*)^2 \otimes (h_1 \wedge h_2)^3, w_1] \\
=& -(r^*)^2 \otimes (h_1 \wedge h_2)^3 \otimes w_1\in\bv_1,
\end{align*}
\item[$\bullet$]
$[\bl_1, \bv_{-1}]$: $$ [r^* \otimes h_1^3, w_1] = -h_1(w_1) r^* \otimes h_1^2\in\bv_0,$$
\item[$\bullet$]
$[\bl_1, \bv_0]$:
\begin{align*}
& [r^*\otimes h_1^3, r^* \otimes h_2^2] \\
=& \langle f, h_1 \wedge h_2\rangle^2 (r^*)^2  \otimes (f^*)^3 \otimes (f \llcorner h_1)\in\bv_1,
\end{align*}
\item[$\bullet$]
$[\bl_{-1}, \bv_0]$: $$[r \otimes w_1^3, r^* \otimes h_1^2] = (h_1(w_1))^2 w_1\in\bv_{-1},$$
\item[$\bullet$]
$[\bl_{-1}, \bv_1]$:
\begin{align*}
& [ r \otimes w_1^3, (r^*)^2 \otimes (f^*)^3 \otimes w_2] \\
=& -\langle w_1\wedge w_2, f^*\rangle r^*\otimes (w_1 \lrcorner f^*)^2\in\bv_0,
\end{align*}
\item[$\bullet$]
$[\bl_{-2}, \bv_1]$: $$[r^2 \otimes f^3, (r^*)^2 \otimes (f^*)^3 \otimes w_3] = w_3\in\bv_{-1}. $$ \end{itemize}
Together with the natural operation of $\bl_0$ on $\bv:= \bv_{-1} \oplus \bv_0 \oplus \bv_1$,
the above operations define a representation of the graded Lie algebra $\bl$ in (2) on $\bv.$
\item[(4)]
The graded Lie algebra $\bl$ is  isomorphic to the graded Lie algebra $\fl$ and the representation $\bv$ of $\bl$ is isomorphic to the irreducible representation $V$ of $\fl$.
\item[(5)] Identify $\fg_0= \fl_0 \oplus \C {\rm Id}_V$ with $\fgl(R) \oplus \fgl(W) = \bl_0 \oplus \C({\rm Id}_W - 3 {\rm Id}_R)$ as Lie algebras by setting $${\rm Id}_W = -3 E -2 {\rm Id}_V, \ \ {\rm Id}_R = -E -{\rm Id}_V \mbox{ and } \fsl(W) = [\fl_0, \fl_0].$$ Then the graded lie algebra $\fg = (\fl \oplus \C{\rm Id}_V) \ltimes V$ is isomorphic to $(\bl \oplus \C {\rm id}_{\bv}) \ltimes \bv$, where ${\rm Id}_{\bv}={\rm Id}_W - 3 {\rm Id}_R$ in the adjoint representation of $(\fgl(R)\oplus\fgl(W))\ltimes\bv$. \end{itemize}
        \end{proposition}

\begin{proof}
(1) can be checked by direct computations, using the natural representations of $\fgl(R)$ and $\fgl(W)$.

(2) can be checked by direct computations, using the formulae given at the beginning of the section.
As an example let us check the Jacobi identity for the triple $(\bl_2, \bl_{-1}, \bl_{-1})$. The Jacobian identity for other triples can be checked similarly.  Take $r, r^*, w_1, w_2, h_1, h_2$ as in the statement of the proposition. Then
\begin{align}\label{e.J1}
& [(r^*)^2\otimes(h_1\wedge h_2)^3, [r\otimes w_1^3, r\otimes w_2^3]] \\
= &[(r^*)^2\otimes(h_1\wedge h_2)^3, r^2\otimes(w_1\wedge w_2)^3] \nonumber\\
=& \langle w_1\wedge w_2, h_1\wedge h_2\rangle^3(-2{\rm Id}_R+{\rm Id}_W). \nonumber
\end{align}
On the other hand,
\begin{align}\label{e.J2}
& [r\otimes w_1^3, [r\otimes w_2^3, (r^*)^2\otimes(h_1\wedge h_2)^3]] \\
=& [r\otimes w_1^3, r^*\otimes(w_2\lrcorner(h_1\wedge h_2))^3] \nonumber\\
=& \langle w_1\wedge w_2, h_1\wedge h_2\rangle^3 {\rm Id}_R \nonumber \\
& -\langle w_1\wedge w_2, h_1\wedge h_2\rangle^2 w_2\lrcorner(h_1\wedge h_2)\otimes w_1. \nonumber
\end{align}
Similarly, we have
\begin{align}\label{e.J3}
& [r\otimes w_2^3, [(r^*)^2\otimes(h_1\wedge h_2)^3, r\otimes w_1^3]] \\
=& [r\otimes w_2^3, -r^*\otimes (w_1\lrcorner(h_1\wedge h_2))^3] \nonumber\\
=& -\langle w_2\wedge w_1, h_1\wedge h_2\rangle^3 {\rm Id}_R \nonumber\\
& + \langle w_2\wedge w_1, h_1\wedge h_2\rangle^2 w_1\lrcorner (h_1\wedge h_2)\otimes w_2. \nonumber
\end{align}
Now we claim that
\begin{align}\label{e.IdW}
&\langle w_1\wedge w_2, h_1\wedge h_2\rangle {\rm Id}_W \\
=& w_2\lrcorner(h_1\wedge h_2)\otimes w_1 -w_1\lrcorner(h_1\wedge h_2)\otimes w_2. \nonumber
\end{align}
The right hand side of \eqref{e.IdW} is an element in $\End(W)$. Its value on $w_1$ is
\begin{align}\label{e.value1}
& (w_2\lrcorner(h_1\wedge h_2)\otimes w_1 -w_1\lrcorner(h_1\wedge h_2)\otimes w_2)(w_1) \\
=& \langle w_1\wedge w_2, h_1\wedge h_2\rangle w_1 - \langle w_1\wedge w_1, h_1\wedge h_2\rangle w_2 \nonumber\\
=& \langle w_1\wedge w_2, h_1\wedge h_2\rangle {\rm Id}_W(w_1), \nonumber
\end{align}
and its value on $w_2$ is
\begin{align}\label{e.value2}
& (w_2\lrcorner(h_1\wedge h_2)\otimes w_1 -w_1\lrcorner(h_1\wedge h_2)\otimes w_2)(w_2) \\
=& \langle w_2\wedge w_2, h_1\wedge h_2\rangle w_1 - \langle w_2\wedge w_1, h_1\wedge h_2\rangle w_2 \nonumber\\
=& \langle w_1\wedge w_2, h_1\wedge h_2\rangle {\rm Id}_W(w_2). \nonumber
\end{align}
If $w_1$ and $w_2$ are linearly dependent, then both hand sides of \eqref{e.IdW} vanish.  If $w_1$ and $w_2$ are linearly independent, then they form a basis of the 2-dimensional vector space $W$ and \eqref{e.IdW} is a consequence of \eqref{e.value1} and \eqref{e.value2}. This proves \eqref{e.IdW}.

Applying \eqref{e.J1}--\eqref{e.IdW}, we can verify the Jacobi identity
\begin{align*}
& [(r^*)^2\otimes(h_1\wedge h_2)^3, [r\otimes w_1^3, r\otimes w_2^3]] \\
 +\, & [r\otimes w_1^3, [r\otimes w_2^3, (r^*)^2\otimes(h_1\wedge h_2)^3]] \\
 +\, & [r\otimes w_2^3, [(r^*)^2\otimes(h_1\wedge h_2)^3, r\otimes w_1^3]] \\
 =\, & 0.
\end{align*}

For (3), note that the definitions of the operations are independent of the choice of the nonzero elements $f$ and $f^*$.
We have to show that for any homogeneous elements $\varphi,\psi\in\bl$ and $u\in\bv$,
\begin{align}\label{e.Jmodule}
[[\varphi,\psi], u]=[\varphi, [\psi, u]]-[\psi, [\varphi, u]].
\end{align} This can be checked by direct computations using the formulae recalled at the beginning of the section.
As an example, let us verify \eqref{e.Jmodule} when $\varphi=r^*\otimes h_1^3\in\bl_1$, $\psi=r\otimes w_1^3\in\bl_{-1}$ and $u=r^*\otimes h_2^2\in\bv_0$. The verification of other cases is similar. In this case, we need to verify
\begin{align}\label{e.Jm}
& [[r^*\otimes h_1^3, r\otimes w_1^3], r^*\otimes h_2^2] \\
=& [r^*\otimes h_1^3, [r\otimes w_1^3, r^*\otimes h_2^2]] -[r\otimes w_1^3, [r^*\otimes h_1^3, r^*\otimes h_2^2]]. \nonumber
\end{align}
The first term of \eqref{e.Jm} is
\begin{align}\label{e.Jm1}
& [[r^*\otimes h_1^3, r\otimes w_1^3], r^*\otimes h_2^2] \\
=& [-h_1(w_1)^3 {\rm Id}_R + h_1(w_1)^2 h_1\otimes w_1, r^*\otimes h_2^2] \nonumber\\
=& h_1(w_1)^3 r^*\otimes h_2^2 - 2h_1(w_1)^2h_2(w_1) r^*\otimes h_1\odot h_2, \nonumber
\end{align}
the second term is
\begin{align}\label{e.Jm2}
& [r^*\otimes h_1^3, [r\otimes w_1^3, r^*\otimes h_2^2]] \\
= & [r^*\otimes h_1^3, h_2(w_1)^2 w_1] \nonumber\\
= & - h_1(w_1) h_2(w_1)^2 r^*\otimes h_1^2, \nonumber
\end{align}
and the third term is
\begin{align}\label{e.Jm3}
& [r\otimes w_1^3, [r^*\otimes h_1^3, r^*\otimes h_2^2]] \\
= & [r\otimes w_1^3, \langle f, h_1\wedge h_2\rangle^2 (r^*)^2\otimes (f^*)^3\otimes f\llcorner h_1] \nonumber\\
= & - \langle f, h_1\wedge h_2\rangle^2 \langle w_1\wedge (f\llcorner h_1), f^*\rangle r^*\otimes (w_1\lrcorner f^*)^2. \nonumber
\end{align}
Applying \eqref{e.Jm1} and \eqref{e.Jm2}, we have
\begin{align}\label{e.Jm4}
& [[r^*\otimes h_1^3, r\otimes w_1^3], r^*\otimes h_2^2] - [r^*\otimes h_1^3, [r\otimes w_1^3, r^*\otimes h_2^2]] \\
=& h_1(w_1) r^*\otimes (h_2(w_1) h_1 - h_1(w_1) h_2)^2 \nonumber\\
=& h_1(w_1) r^*\otimes (w_1\lrcorner (h_1\wedge h_2))^2. \nonumber
\end{align}
If $h_1$ and $h_2$ are linearly dependent, then \eqref{e.Jm} is a direct consequence of \eqref{e.Jm3} and \eqref{e.Jm4}.
Assume that $h_1$ and $h_2$ are linearly independent. Then
\begin{align}\label{e.scaler}
f^*=\lambda h_1\wedge h_2 \mbox{ for some } \lambda\in\C^*.
\end{align}
It follows that
\begin{align}\label{e.coeffi}
& \langle w_1\wedge (f\llcorner h_1), f^*\rangle \\
=& \lambda\langle w_1\wedge(f\llcorner h_1), h_1\wedge h_2\rangle \nonumber\\
=& -\lambda\langle f\llcorner h_1, w_1\lrcorner(h_1\wedge h_2)\rangle  \nonumber\\
=& -\lambda(f\llcorner h_1)\llcorner(w_1\lrcorner(h_1\wedge h_2)) \nonumber\\
=& -\lambda\langle f, h_1\wedge(w_1\lrcorner(h_1\wedge h_2)) \nonumber\\
=& -\lambda\langle f, h_1\wedge (h_2(w_1)h_1-h_1(w_1)h_2) \nonumber\\
=& \lambda h_1(w_1)\langle f, h_1\wedge h_2\rangle \nonumber\\
=& h_1(w_1). \nonumber
\end{align}
Applying \eqref{e.scaler} and \eqref{e.coeffi} into \eqref{e.Jm3}, we have
\begin{align}\label{e.Jm5}
& [r\otimes w_1^3, [r^*\otimes h_1^3, r^*\otimes h_2^2]] \\
= & -\lambda^{-2} h_1(w_1) r^*\otimes (w_1\lrcorner (\lambda h_1\wedge h_2))^2 \nonumber\\
= & -h_1(w_1) r^*\otimes (w_1\lrcorner(h_1\wedge h_2))^2. \nonumber
\end{align}
Then \eqref{e.Jm} follows from \eqref{e.Jm4} and \eqref{e.Jm5}.

Now let us verify (4). It is easy to check that $\bl_{-2} \oplus \bl_{-1} \oplus \bl_0$ is isomorphic to $\fl_{-2} \oplus \fl_{-1} \oplus \fl_0$ as graded Lie algebras.  Then the isomorphism $\fl \cong \bl$ follows from the definition of the universal prolongation and Proposition \ref{p.kim}, because $\dim \fl = \dim \bl$.
Since $\dim V = \dim \bv =7$ and there is only one irreducible representation of $\fl$ with dimension less than or equal to $7$ up to isomorphism, the isomorphism $V \cong \bv$ as  representations of $\fl \cong \bl$ is immediate.

(5) is straight-forward once we have (4).
\end{proof}

We note the following, whose proof is easy.

\begin{lemma}\label{l.Pspan}
In Proposition \ref{p.repre} (i), let $W^0 \subset W$ be any nonempty Euclidean open subset in $W$.
 Then the kernel of the Lie bracket $\wedge^2 \fl_{-1} \to \fl_{-2}$ is the linear span of
$$\bigcup_{w \in W^o}  (R \otimes w^3) \wedge (R \otimes (w^2 \odot W)) $$
in $\wedge^2(R \otimes \Sym^3 W) = \wedge^2 \fl_{-1}.$  \end{lemma}

\section{Geometry of the surface $\BS \subset \BP^5$}\label{s.S}

In this section, we study the geometry of some rational ruled surfaces, including the projective  surface  $\BS \subset \BP^5$ discussed in Section \ref{s.I}. These surfaces are  related to
the VMRT of the ${\rm G}_2$-horospherical variety $\bX$ at a general point, as we see  in Lemma \ref{l.deform} and Corollary \ref{c.bX}.

\begin{definition}\label{d.scroll}
For two natural numbers $a \leq b$, the sections of the vector bundle $\sO(a) \oplus \sO(b)$ on $\BP^1$ induces an embedding of the
 ruled surface $\BP (\sO(-a) \oplus \sO(-b))$ into
 $ \BP^{a+b+1}$ the  image of which we denote by $S(a,b)$. This projective surface
 $S(a,b)$ is called the {\em rational normal scroll} of type $(a,b)$ (see Section 1 of \cite{EH}).
 Denote by $L_{ab}$ the line bundle on $S(a,b)$ defining the embedding into $\BP^{a+b+1}$.
 \end{definition}

 We recall the following properties of rational normal scrolls. (i), (ii) and (iii)  are elementary and  can be seen easily from Section 1 of \cite{EH}  and Section V.2 of \cite{Ha}. (iv) is from Lemma 2.1 of \cite{EH}.

 \begin{lemma}\label{l.ruled}
\begin{itemize} \item[(i)]  The surface $S(a,b)$ is a rational ruled surface, equipped with a $\BP^1$-fibration $\phi: S(a,b) \to \BP^1$.
\item[(ii)] There exists a  section $E \subset S(a,b)$ of $\phi$
with the self-intersection number $E\cdot E = a-b \leq 0.$
\item[(iii)] The degree of $S(a,b) \subset \BP^{a+b+1}$ is $L_{ab}\cdot L_{ab} = a+b$
and when $F \subset S(a,b)$ is a fiber of $\phi$, we have the numerical equivalence $L_{ab} \equiv E + b F$.
\item[(iv)] The surface $S(a,b) \subset \BP^{a+b+1} $ is defined by quadratic equations. \end{itemize} \end{lemma}

We are mainly interested in  $S(1,3)$ and $S(2,2)$. First, we have the following description of polarized deformations of these surfaces.

\begin{proposition}\label{p.deform}
Let $\pi:\mathcal{S}\to\Delta$ be a smooth family of projective surfaces and $\mathcal{L}$ be a $\pi$-ample  line bundle on $\mathcal{S}$. Define $\sS_t:= \pi^{-1}(t)$ and $\sL_t := \sL|_{\sS_t}$ for each $t \in \Delta$. 
\begin{itemize}
\item[(i)] If $(\mathcal{S}_t, \mathcal{L}_t)$ is isomorphic to  $(S(2,2), L_{22})$ for $t\neq 0$, then $(\mathcal{S}_0, \mathcal{L}_0)$ is isomorphic to either $(S(2,2), L_{22})$ or $(S(1,3), L_{13})$.
\item[(ii)] If $(\mathcal{S}_t, \mathcal{L}_t)$ is isomorphic to  $(S(1,3), L_{13})$ for $t\neq 0$, then so is $(\mathcal{S}_0, \mathcal{L}_0)$.
\end{itemize}
\end{proposition}
\begin{proof}
By the classification of surfaces,  the central fiber $ \mathcal{S}_0$ in both (i) and (ii) is  a rational ruled surface  with a $\BP^1$-fibration $\phi: \mathcal{S}_0 \to \BP^1,$ which is the limit of the $\BP^1$-fibration on $\mathcal{S}_t$. We have a section $E \subset \mathcal{S}_0$ of $\phi$ with $e:= - E \cdot E \geq 0$ as in Lemma \ref{l.ruled}.     The line bundle  $\mathcal{L}_0$ has degree 1 on a fiber $F$ of $\varphi$. Thus we have the numerical equivalence $\mathcal{L}_0 \equiv E + c F$ for some integer $c$.  As $\mathcal{L}_0\cdot \mathcal{L}_0 = \mathcal{L}_t \cdot
\mathcal{L}_t = 4$ by Lemma \ref{l.ruled}, we have $$4 = (E + cF) \cdot (E + cF) = -e + 2c.$$
From $1 \leq E \cdot \mathcal{L}_0 = -e + c = 4-c $, we obtain $c =2, e=0$ or $c =3, e=2$.
 From the well-known classification of rational ruled surfaces,   the integers $e$ and $c$ determine the isomorphism class of $(\mathcal{S}_0, \mathcal{L}_0).$ This proves (i). For (ii), it suffices to exclude the possibility that $\mathcal{S}_0 \cong S(2,2)$. But this follows from the infinitesimal deformation rigidity of $S(2,2) \cong \BP^1 \times \BP^1$ (for example, Korollar 6.4 of \cite{Br}). \end{proof}

For details on the following notion, Section 12.1 of \cite{IL} is a good reference.

\begin{notation}\label{n.ff}
Let $U$ be a vector space and let $Y \subset \BP U$ be a (locally closed) complex
submanifold. For each point $y \in Y$, we have the $k$-th osculating space
$\sT^k_y \subset U$ and the $k$-th normal space $\sN^k_y$ of $Y$ at $y$
(see Section 12.1 of \cite{IL}) such that
\begin{itemize}
\item[(1)] $\sT^0_y\subset U$ is the 1-dimensional space corresponding to $y \in \BP U$;
    \item[(2)] $\sT^1_y \subset U$ is the affine tangent space of $Y $ at $y$;
    \item[(3)] $\sN^k_y = \sT^k_y/\sT^{k-1}_y$; and
    \item[(4)] if $\dim \sN^i_z$ is constant for all $z \in Y$ in a neighborhood of $y$ for all $i \leq k$, then we have a surjective homomorphism
        $${\rm FF}^k_y: \Sym^k T_y Y \to N^{(k)}_y:= \Hom(\sT^0_y, \sN^k_y),$$
        called the  $k$-th fundamental form. For convenience, we write ${\rm II}_y$ for ${\rm FF}^2_y$  and ${\rm III}_y$ for ${\rm FF}^3_y$.
        \end{itemize} \end{notation}

\begin{lemma}\label{l.scroll}
Let $S \subset \BP^5$ be either $S(1,3)$ or $S(2,2)$. At a general point $y \in S$,  the system of quadrics $| \II_y |\subset \BP  (\Sym^2 T^*_{y} S)$ defined by the second fundamental form has a single base point on $\BP T_y S$. Consequently, we have $\dim N_{y}^{(2)} =2$. \end{lemma}

\begin{proof}
 As $S$ is defined by
qudratic equations in $\BP^5 $, lines of $\BP^5$ through $y$ that have contact order 2 with $S$ must be contained in $S$. This implies that  the base locus of the system $| \II_{y}|$ must be tangent vectors to lines on $S$ through $y$.
It is easy to check that  there exists a unique line on $S$ passing through a general point, namely, the
 fiber of the ruled surface $S(a,b) \to \BP^1$.  Thus the base locus of  the system of quadrics $| \II_{y}|$ on $\mathbb{P}T_{y} S$ is a single point. It follows that  $\dim | \II_{y}| =1$ and $\dim |\II_y| + 1= \dim N_{y}^{(2)} =2.$ \end{proof}

 We have the following characterization of surfaces of minimal degree, due to del Pezzo (Theorem 1 in \cite{EH}).

 \begin{theorem}\label{t.DelPezzo}
 Let $S \subset \BP^{5}$ be a linearly nondegenerate irreducible projective surface of degree  $4$. Then it is projectively equivalent to one of the following:
 \begin{itemize}
 \item[(1)] a cone over a rational normal curve of degree $4$,
 \item[(2)] the Veronese surface $v_2(\BP^2) \subset \BP^5$,
 \item[(3)] rational normal scrolls $S(1,3)$ or $S(2,2)$.
\end{itemize}
\end{theorem}

We have the following two corollaries of Theorem \ref{t.DelPezzo}.

\begin{corollary}\label{c.DelPezzo}
Let $S \subset \BP^5$ be an irreducible surface of degree at most $4$. Assume that  both the second
and the third fundamental forms of $S$ at a general point of $S$ are isomorphic to those of $S(1,3)$.
Then $S$ is linearly nondegenerate in $\BP^5$, and it is projectively isomorphic to $S(1,3)$ or $S(2,2)$.
\end{corollary}
\begin{proof}  Let $y$ be a general point of $S(1,3)$. Then Lemma \ref{l.scroll} shows that $\dim N_{y}^{(2)} =2$  and consequently $\III_{y}\neq 0$. By the assumption on fundamental forms, $S$ is linearly nondegenerate in $\BP^5$. Then by Proposition 0 in  \cite{EH}, the degree of $S\subset\BP^5$ is at least 4, hence equal to 4.
So $S$ must be one of
the surfaces (1)-(3) in Theorem \ref{t.DelPezzo}. The second fundamental
form of the cone (1) cannot be equal to that of a nondegenerate nonsingular surface (e.g. Proposition 4.3.1 and Theorem 4.3.6 of \cite{IL}). The
third fundamental form of the Veronese surface (2) vanishes (e.g.  Section 12.2 of \cite{IL}). Thus the only possibility is (3). \end{proof}

\begin{corollary}\label{c.deform}
Let $\{S_t \subset \BP^5, \ t \in \Delta\}$ be a well-defined family of irreducible reduced algebraic cycles of dimension 2 and degree 4 in the sense of Definition I.3.10 of \cite{Ko}. Assume that  the fiber $S_0$ at $0 \in \Delta$ is linearly nondegenerate, and it is not (1) in Theorem \ref{t.DelPezzo}.
\begin{itemize} \item[(i)] If $S_t$ is projectively isomorphic to $S(2,2)$ for $t \neq 0$, then $S_0$ is projectively isomorphic to
either $S(2,2)$ or $S(1,3)$.
\item[(ii)] If $S_t$ is
projectively isomorphic to $S(1,3)$ for $t \neq 0,$ then so is $S_0$.
\end{itemize} \end{corollary}

\begin{proof}
By Theorem \ref{t.DelPezzo}, the surface $S_0$ is (2) or (3). Then the family of cycles $S_t \subset \BP^5$ is a flat family of nonsingular surfaces
by Theorem I.6.5 of \cite{Ko}. Thus the result follows from Proposition \ref{p.deform}.   \end{proof}

Now we relate the surface $S(1,3)$ to the graded Lie algebra in Proposition \ref{p.repre}. The following is easy to check.

\begin{proposition}\label{p.BS}
Let $R$ be a 1-dimensional vector space and let $W$ be a 2-dimensional vector space.
Consider the affine variety  $\widehat{\BS} \subset W \oplus (R \otimes \Sym^3 W)$ defined by
$$\widehat{\BS} := \{ c \cdot w + r \otimes w^3, \ c \in \C, r \in R, w \in W\}.$$ The corresponding projective variety $\BS \subset \BP (W \oplus (R \otimes \Sym^3 W))$ is biregular to $S(1,3)$. \end{proposition}

Recall the following from Satz 1.10 of \cite{Br}.

\begin{proposition}\label{p.Br}
The automorphism group of the surface $S(1,3) \cong \BS$ is connected and has dimension $7$. \end{proposition}

\begin{proposition}\label{p.autS}
Let $\bG_0 \subset {\rm GL}(W \oplus (R \otimes \Sym^3 W))$ be the image of the effective representation  of  $$({\rm GL}(R) \times {\rm GL}(W)) \ltimes (R^* \otimes \Sym^2 W^*) $$ on $W \oplus (R \otimes \Sym^3 W).$  \begin{itemize}
\item[(1)] The linear  group  $\bG_0$   is the connected component of the linear automorphism group of the affine variety $\widehat{\BS} \subset  W \oplus (R \otimes \Sym^3 W).$
    \item[(2)] The induced action of $\bG_0 $ on $\BS$  has two orbits, the closed orbit $\BP W \subset \BS$ and its complement.  \end{itemize} \end{proposition}

\begin{proof}
It is easy to see that $\bG_0$ preserves $\widehat{\BS}.$ By Proposition \ref{p.Br},
the linear automorphism group of the affine variety $\widehat{\BS}$ has dimension at most $8 = \dim \bG_0$. Thus its connected component of identity coincides with $\bG_0$, proving (1). (2) is easy to check. \end{proof}

\begin{lemma}\label{l.IIz}
At a point $z \in \BP W \subset \BS$ of the closed orbit in Proposition \ref{p.autS} (2), the base loci of the system of quadrics $| \II_y|$ on $\BP T_y S$ has two distinct points. Consequently, it is not isomorphic to $|\II_y|$ at a point $y \in \BS \setminus \BP W$. \end{lemma}

\begin{proof}
As in the proof of Lemma \ref{l.scroll}, this follows from the fact that there are two distinct lines on $\BS$ through $z$: the fiber of $S(1,3) \to \BP^1$ and the line $\BP W \subset \BS$. \end{proof}

The following can be checked by a direct computation.

\begin{lemma}\label{l.ST}  In Proposition \ref{p.BS},
fix a point $s := w + r \otimes w^3 \in \widehat{\BS}$ with  $r \neq 0 \neq w$ and set $W_{\flat} := \C w \subset W.$ Let $\sT^i_s$ be the $i$-th osculating space of $\BS$ at the point $[s] \in \BS$ corresponding to $s$ in Notation \ref{n.ff}.
Then we have the following identification \begin{eqnarray*}
\sT^0_s & = & \C s \\
\sT^1_s &=& \{v+3r\otimes w^2\odot v+\lambda r\otimes w^3, \lambda\in\C, v\in W\} \\
\sT^2_s&=& W  \oplus (R \otimes (W_{\flat} \odot \Sym^2 W)) \\
\sT^3_s & =& W \oplus (R \otimes \Sym^3 W) \end{eqnarray*}
as subspaces of $\fg_{-1} = W \oplus (R \otimes \Sym^3 W)$ in Proposition \ref{p.repre}.
In particular, the vector space $\sT^1_s$ satisfies
$$W_{\flat}\oplus(R\otimes\Sym^3W_{\flat}))\subsetneqq \sT^1_s\subsetneqq W\oplus(R\otimes\Sym^2W_{\flat}\odot W),$$
and the intersections of subspaces of $\fg_{-1}$ satisfy
$$W \cap \sT^1_s = W_{\flat} \mbox{ and } W \cap \sT^2_s = W.$$ \end{lemma}

The next proposition is immediate from Lemma \ref{l.ST}.

\begin{proposition}\label{p.ST}
Let us use the identification $\fg_{-1} = W \oplus (R \otimes \Sym^3 W)$ from Proposition \ref{p.repre} and regard $\widehat{\BS}$ as a subvariety in $\fg_{-1}$.
In Lemma \ref{l.ST}, write $W_{\sharp} := W/W_{\flat}$.
 The exact sequence of $\bG_0$-modules
\begin{equation}\label{e.G0} 0 \to W \to \fg_{-1} \to R \otimes \Sym^3 W \to 0 \end{equation}
induces an isomorphism $\sT^0_s \cong R \otimes \Sym^3 W_{\flat}$
and short exact sequences
$$0 \to W_{\flat} \to \sT^1_s \to R \otimes (\Sym^2 W_{\flat} \odot W) \to 0, $$
$$0 \to W \to \sT^2_s \to R \otimes (W_{\flat} \odot \Sym^2 W) \to 0,$$ $$ 0 \to W \to \sT^3_s \to R \otimes \Sym^3 W \to 0.$$
Consequently, (\ref{e.G0}) induces   the following exact sequences for the $i$-th normal space $\sN^{i}_s = \sT^{i+1}_s/\sT^i_s$ of $\BS$ at $s$,
$$ 0 \to W_{\flat} \to \sN^1_s \to R \otimes (\Sym^2 W_{\flat} \odot W_{\sharp}) \to 0,$$
$$0 \to W_{\sharp} \to \sN^2_s \to R \otimes (W_{\flat} \odot \Sym^2 W_{\sharp}) \to 0,$$
as well as an isomorphism $\sN^3_s \cong R \otimes \Sym^3 W_{\sharp}.$ \end{proposition}

\begin{proposition}\label{p.span}
Using the identification $\fg_{-1} = W \oplus (R \otimes \Sym^3 W)$ from Proposition \ref{p.repre} and  $\widehat{\BS}\subset \fg_{-1},$
let $\Xi \subset \wedge^2 \fg_{-1}$ be the subspace spanned by subspaces of the form $\wedge^2 P$ where  $P \subset \fg_{-1}$ is any 2-dimensional subspace tangent to $\widehat{\BS}$. Then $\Xi = {\rm Ker}(\psi)$ where $\psi : \wedge^2 \fg_{-1} \to \fg_{-2}$ is the Lie bracket.  \end{proposition}

\begin{proof}
We skip the proof of the inclusion $\Xi \subset {\rm Ker}(\psi)$, which can be checked either by a direct computation or using Proposition \ref{p.detP} below combined with Corollary \ref{c.bX}.

To prove the inclusion ${\rm Ker}(\psi) \subset \Xi,$ we identify  $R$ with $\C$ and view $\wedge^2 \fg_{-1} = \wedge^2 (W \oplus \Sym^3 W)$ as an $\fsl(W)$-module.
Both $\Xi$ and ${\rm Ker}(\psi)$ are preserved under the $\fsl(W)$-action.
We know that ${\rm Ker}(\psi)$ has codimension 1 in $\wedge^2 \fg_{-1}.$
In the $\fsl(W)$-module decomposition
\begin{equation}\label{e.decompo} \wedge^2 \fg_{-1} = \wedge^2 W \oplus (W \wedge \Sym^3 W) \oplus \wedge^2 (\Sym^3 W),\end{equation} the last factor is decomposed as $$\wedge^2 (\Sym^3 W) \  \cong \  \Sym^4 W \oplus \C$$ as $\fsl(W)$-modules.
 Since $ W \subset \widehat{\BS}$
  and $ \{  w^3, w \in W\}  \subset \widehat{\BS},$  it is easy to see that $\wedge^2 W$ and the $\Sym^4 W$-factor of $\wedge^2 (\Sym^3 W)$ in (\ref{e.decompo})  are contained in $\Xi$. It remains to show that $\Xi$ contains $ W \wedge \Sym^3 W$
  modulo elements of $\wedge^2 W$ and the $\Sym^4 W$-factor of $\wedge^2 (\Sym^3 W)$.

  To check this,   regard $\Sym^3 W$ as a subspace of $W \otimes W \otimes W$.
  For any $v, w \in W$, consider the following arc on $\widehat{\BS}$,
  $$s_t := (v + t w) +  (v + tw) \otimes (v + tw) \otimes (v + tw), \ t \in \Delta.$$
  Then $\Xi$ contains $$s_0 \wedge \frac{{\rm d} s_t}{{\rm d} t}|_{t=0} = (v +  v \otimes v \otimes v) \wedge (w + \square),$$ where $\square =v \otimes v \otimes w + v \otimes w \otimes v +  w \otimes v \otimes v.$ It follows that $\Xi$ contains, modulo elements of $\wedge^2 W$ and the $\Sym^4 W$-factor of $\wedge^2 (\Sym^3 W),$ $$ ( v \otimes v \otimes v)  \wedge w + v \wedge \square,$$ which can be expanded in $\wedge^2(W \otimes W^{\otimes 3}) \subset (W \oplus W^{\otimes 3}) \otimes (W \otimes W^{\otimes 3})$ as    \begin{eqnarray*}  v \otimes v \otimes v \otimes w -  w \otimes v \otimes v \otimes v + v \otimes v \otimes v \otimes w +  v \otimes v \otimes w \otimes v  \\  +  v \otimes w \otimes v \otimes v  - v \otimes v \otimes w \otimes v -  v \otimes w \otimes v \otimes v - w \otimes v \otimes v \otimes v. \end{eqnarray*}  This is $2 ( v  \otimes v \otimes v) \wedge w.$ Since elements of the form $v \otimes v\otimes v$ span $\Sym^3 W$, we see that $W \wedge \Sym^3 W$ lie in $\Xi$. \end{proof}

One consequence of Proposition \ref{p.span} is the following.

\begin{proposition}\label{p.span2} Let $V$ be a vector space and let
$S \subset \BP V$ be any irreducible subvariety with the affine cone $\widehat{S} \subset V$. Define $\Xi_S \subset \wedge^2 V$ as the linear subspace spanned by $$\{ \wedge^2 P \mid \mbox{ 2-dimensional subspace } \dim  P  \subset V \mbox{ tangent to }  \widehat{S} \}.$$
   \begin{itemize} \item[(i)] For the scroll $S(2,2) \subset  \BP H^0(S(2,2), L_{22})^* ,$ the subspace $\Xi_{S(2,2)}$  is a hyperplane in $\wedge^2 H^0(S(2,2), L_{22})^*$  and
 the projection $$\wedge^2 H^0(S(2,2), L_{22})^* \to \C$$ modulo $\Xi_{S(2,2)}$ defines a symplectic form on $H^0(S(2,2), L_{22})$.
 \item[(ii)] Let $S \subset \BP V$ be a surface which is the image of $S(2,2)$ under a  linear projection $H^0(S(2,2), L_{22})^* \to V$  with nonzero kernel. Then $\Xi_S = \wedge^2 V$.
      \item[(iii)] For the scroll $S(1,3) \subset \BP H^0(S(1,3), L_{13})^*$, the subspace $\Xi_{S(1,3)}$  is a hyperplane in $\wedge^2 H^0(S(1,3), L_{13})^*$  and
 the projection $$\wedge^2 H^0(S(1,3), L_{13})^* \to \C$$ modulo $\Xi_{S(1,3)}$ defines a antisymmetric form on $H^0(S(1,3), L_{13})$, the null space of which is exactly the two-dimensional subspace corresponding to $W$ in Proposition \ref{p.BS}.
 \item[(iv)]  Let $S \subset \BP V$ be a surface which is  the image of $S(1,3)$ under a linear projection $H^0(S(1,3), L_{13})^* \to V$  with nonzero kernel such that $S(1,3)$ is the normalization of $S$ by this projection. Then $\Xi_S = \wedge^2 V$.
 \end{itemize}
    \end{proposition}

    \begin{proof} (i) is from Proposition 4 of \cite{Hw97}. 

    For (ii), let $V' \subset H^0(S(2,2), L_{22})^* $  be the kernel of the projection to $V$. Then $V' \wedge H^0(S(2,2), L_{22})^*$  is the kernel of the induced projection $\wedge^2 H^0(S(2,2), L_{22})^* \to \wedge^2 V$. To prove (ii), it suffices to show that $\Xi_{S(2,2)}$ surjects to $\wedge^2 V$, in other words, the kernel $V' \wedge H^0(S(2,2), L_{22})^*$ is not contained in the hyperplane $\Xi_{S(2,2)}$. But if $V' \wedge H^0(S(2,2), L_{22})^* \subset \Xi_{S(2,2)}$, then $V'$ is in the null space of  the symplectic form in (i), a contradiction.

(iii) is from Proposition \ref{p.span} and the property of the Lie bracket $\wedge^2 \fg_{-1} \to \fg_{-2}$ in Definition \ref{d.fg}.

For (iv), let $V' \subset H^0(S(1,3), L_{13})^* $  be the kernel of the projection to $V$. Then $V' \wedge H^0(S(1,3), L_{13})^*$  is the kernel of the induced projection $\wedge^2 H^0(S(2,2), L_{22})^* \to \wedge^2 V$. To prove (iv), it suffices to show that $\Xi_{S(2,2)}$ surjects to $\wedge^2 V$, in other words, the kernel $V' \wedge H^0(S(1,3), L_{13})^*$ is not contained in the hyperplane $\Xi_{S(1,3)}$. Suppose $V' \wedge H^0(S(1,3), L_{13})^* \subset \Xi_{S(1,3)}$, then $V'$ is in the null space $W$ of  the antisymmetric form in (iii). From the description of $S(1,3)$ in Proposition \ref{p.BS}, the projection of $S(1,3)$ from a subspace in $W$ must contract a line in $S(1,3)$. So  $S(1,3)$ cannot be the normalization of $S$ by this projection, a contradiction.  \end{proof}

\section{Vanishing of sections of vector bundles}\label{s.vanishing}
In this section, we prove the vanishing of sections of certain vector bundles on a manifold satisfying certain assumptions. This vanishing result is needed  to apply Theorem \ref{t.HL} in Section \ref{s.proof}. Throughout this section, we work in the following setting.

\begin{setup}\label{setup}
 Let $\fg$ be as in Definition \ref{d.fg} and let $\bG_0$ be as in Proposition \ref{p.autS}.   We consider a complex manifold $M$ of dimension 7 with  a vector subbundle $\sD \subset TM$ of rank 6 and a closed submanifold $\sC \subset \BP \sD$ such that
\begin{itemize}
\item[(i)]
the filtration $$ F^{\bullet} = (F^{-1}= \sD  \subset F^{-2}= TM)$$ on $M$ is a filtration of type $\fg_-= \fg_{-1} \oplus \fg_{-2}$;
\item[(ii)] there is a $\bG_0$-structure $\sA \subset {\rm grFr}(F^{\bullet})$ subordinate to the filtration $F^{\bullet}$; and
    \item[(iii)]  the isomorphism from $\fg_{-1}$ to $\sD_x$ induced by any  element of $\sA_x$ at any $x \in M$ sends $\BS \subset \BP \fg_{-1}$ isomorphically to the fiber  $\sC_x \subset \BP \sD_x$ of the projection $\sC \to M$ at $x \in M$. \end{itemize}
Define the following vector bundles associated to $\sA$ by the natural $\bG_0$-action:
$$\sV_k := \sA \times^{\bG_0} V_k \mbox{ and } \sG_k := \sA \times^{\bG_0} \fg_k$$ for $-2 \leq k \leq 2$ with $V_{-2} = V_2 =0.$
Let $\sL_k$ be the quotient bundle $\sG_k/\sV_k$.
Then writing $\sW = \sV_{-1}$ and $\sQ = \sL_{-2}$, we can find a line bundle $\sR$ on $M$ satisfying the following relations analogous to Proposition \ref{p.repre}:
\begin{eqnarray*}
\sQ & = & \sR^{\otimes 2} \otimes (\wedge^2 \sW)^{\otimes 3} \\
\sL_{-2} & = & \sQ \\
\sL_{-1} & = & \sR \otimes \Sym^3 \sW\\
\sL_0 & = & \sO_M \oplus (\sW^* \otimes \sW) \\
\sL_1 & = & \sR^* \otimes \Sym^3 \sW^* \\
\sL_2 & = & \sQ^*\\
\sV_{-1} & =& \sW \\
\sV_0 & = & \sR^* \otimes \Sym^2 \sW^* \\
\sV_1 &=& \sQ^* \otimes \sW.\end{eqnarray*}
\end{setup}

\begin{remark}
In Setup \ref{setup}, each $\sL_i$ corresponds to $\fl_i$ or  $\bl_i$ in Proposition \ref{p.repre}, except $\sL_0$ which corresponds to $\bl_0 \oplus \C {\rm Id}_V$. This is to simplify our notation below. \end{remark}

We make the following additional assumption.

\begin{assumption}\label{assumption}
Fix a base point $o \in \BP^1$.
In Setup \ref{setup}, there exist
\begin{itemize}
\item[(1)] a nonempty Zariski-open subset $M' \subset M$;
\item[(2)]  a nonempty Euclidean open subset $\sW_x^o \subset \sW_x \setminus \{0\} $ for each $x \in M'$; and
    \item[(3)] a holomorphic map  $f_a: \BP^1 \to M$ for each $a \in \sW^o_x, x \in M'$, such that
\begin{itemize} \item[(i)] $f_a(o) = x$;  \item[(ii)] there are isomorphisms
$$f_a^* \sQ \cong \sO(1) \mbox{ and } f_a^*\sR \cong \sO(-1)$$ of line bundles on $\BP^1$; and  \item[(iii)] $f_a^* \sW$ contains a unique line subbundle $f^1_a \sW$ of degree 1 on $\BP^1$ such that $f_a^*\sW \cong f^1_a \sW \oplus \sO_{\BP^1}$ and the fiber of $f^1_a\sW$ at $o$ corresponds to the subspace $\C a$ of $\sW_x$.
 \end{itemize} \end{itemize}
 \end{assumption}

We consider the following conditions for vector bundles on $M$.

\begin{definition}\label{d.generation}
In Assumption \ref{assumption},   for a vector bundle $\sU$ on $M$ and an integer $k$, denote by $f^k_a \sU$ the unique subbundle of $f_a^* \sU$ spanned by line subbundles of degree $\geq k$ on $\BP^1$. We define the following notions for a positive integer $m.$
\begin{itemize}
    \item[(1)] We say that $\sU$ is \emph{ $m$-tuply  $\sO(k)$-generated} if for a general point $x \in M',$
    $$\bigcup_{a_1, \ldots, a_m} ( \bigcap_{1 \leq i \leq m} (f^k_{a_i} \sU)_o) \mbox{ spans the fiber } f_a^* \sU_x,$$
    where the union is taken over all pairwise linearly independent $m$-tuples $a_1, \ldots, a_m$ of points in $\sW_x^o$. We say that $\sU$ is $\sO(k)$-generated if it is $m$-tuply $\sO(k)$-generated for $m=1$.
    \item[(3)] We say that $\sU$ is \emph{$m$-tuply $\sO(k)$-null} if for a general point $x \in M',$
    $$\bigcap_{1 \leq i \leq m} (f^k_{a_i} \sU)_o = 0$$ for any pairwise linearly independent $m$-tuples $a_1, \ldots, a_m$ of points in $\sW_s^o$. We say that $\sU$ is $\sO(k)$-null if it is $m$-tuply $\sO(k)$-null for $m=1$. \end{itemize}
    \end{definition}

The following three lemmata are easy to see.

\begin{lemma}\label{l.compare}
In Definition \ref{d.generation}, let $\sU_1$ (resp. $\sU_2$) be a vector bundle on $M$ which is $m$-tuply $\sO(k)$-generated (resp. $m$-tuply $\sO(k)$-null). Then the following holds. \begin{itemize}
\item[(1)] $\sU_1$ is $m'$-tuply $\sO(k')$-generated for all $1 \leq m' \leq m$ and $k' \leq k$.
    \item[(2)] $\sU_2$ is $m'$-tuply $\sO(k')$-null for all $m' \geq m \geq 1$ and $k' \geq k$.
        \item[(3)] $H^0(M, \sU_1^* \otimes \sU_2) = {\rm Hom}(\sU_1, \sU_2) =0.$ \end{itemize} \end{lemma}

\begin{lemma}\label{l.line}
In Definition \ref{d.generation}, let $\sU$ be a line bundle on $M$ and let $k$ be an integer. Then the following three conditions are equivalent.
\begin{itemize}
\item[(i)] $\sU$ is $m$-tuply $\sO(k)$-generated for some positive integer $m$.
\item[(ii)] $\sU$ is $m$-tuply $\sO(k)$-generated for any positive integer $m$.
\item[(iii)] $\deg f_a^* \sU \geq k$ for all $a \in \sW^o_x$. \end{itemize}
\end{lemma}

\begin{lemma}\label{c.split}
In Assumption \ref{assumption},  for each $x \in M'$ and $a \in \sW^o_x$,  we have the following isomorphisms of vector bundles on $\BP^1.$
\begin{eqnarray*}
f_a^*\sL_{-2} & \cong & \sO(1) \\
f_a^*\sL_{-1} & \cong & \sO(2) \oplus \sO(1) \oplus \sO \oplus \sO(-1) \\
f_a^*\sL_0 & \cong & \sO(1) \oplus \sO^{\oplus 3} \oplus \sO(-1) \\
f_a^* \sL_1 & \cong & \sO(1) \oplus \sO \oplus \sO(-1) \oplus \sO(-2) \\
f_a^* \sL_2 & \cong & \sO(-1) \\
f_a^* \sV_{-1} & \cong & \sO(1) \oplus \sO\\
f_a^* \sV_0 & \cong & \sO(1) \oplus \sO \oplus \sO(-1) \\
f_a^* \sV_1 & \cong & \sO \oplus \sO(-1). \end{eqnarray*}
\end{lemma}

\begin{proposition}\label{p.null} In Assumption \ref{assumption},
we have the following properties of vector bundles on $M$.
\begin{itemize}
\item[(i)] $\sL_{-2}$ is $\sO(2)$-null.
\item[(ii)] $\sL_{-1}$ is $\sO(3)$-null and 4-tuply $\sO$-null.
\item[(iii)] $\sL_0$ is $\sO(2)$-null and 2-tuply $\sO(1)$-null.
\item[(iv)] $\sL_1$ is $\sO(2)$-null and 2-tuply $\sO(1)$-null.
\item[(v)] $\sL_2$ is $\sO$-null.
\item[(vi)] $\sV_{-1}$ is $\sO(2)$-null and 2-tuply $\sO(1)$-null.
\item[(vii)] $\sV_0$ is $\sO(2)$-null and 2-tuply $\sO(1)$-null.
\item[(viii)] $\sV_1$ is $\sO(1)$-null.
\end{itemize} \end{proposition}

\begin{proof}
The assertions in Proposition \ref{p.null} that a vector bundle is  $m$-tuply $\sO(k)$-null with $m=1$ are obvious from Lemma \ref{c.split}.
Let us prove the assertions with $m\geq 2$.

For any $x \in M'$ and $a \in \sW^o_x$, we know from Assumption \ref{assumption} (3) that $$(f^0_a \sL_{-1})_x = \sR_x \otimes (a \odot \Sym^2 \sW_x)$$ under the identification $\sL_{-1} = \sR \otimes \Sym^3 \sW$. Then for any pairwise linearly independent $3$-tuple $a_1, a_2, a_3 \in\sW_x^o$, $$(f^0_{a_1} \sL_{-1})_x \cap (f^0_{a_2} \sL_{-1})_x \cap (f^o_{a_3} \sL_{-1})_x = \sR_x \otimes (a_1 \odot a_2 \odot a_3).$$ Thus for any pairwise linearly independent 4-tuple $a_1, a_2, a_3, a_4 \in \sW_x^o$, we have $$\bigcap_{1 \leq i \leq 4} (f^0_{a_i} \sL_{-1})_x = 0,$$ which proves (ii).

For any $x \in M'$ and $a \in \sW^o_x$, we have $$(f^1_a \sL_0)_x \subset \sW^*_x \otimes a$$ under the identification $\sL_0 \cong (\sW^* \otimes \sW) \oplus \sO_M.$ Thus for any linearly independent $a_1, a_2 \in \sW^o_x$, we have
$$(f^1_{a_1} \sL_0)_x \cap (f^1_{a_2} \sL_0)_x =0,$$ which proves (iii).

For any $x \in M'$ and $a \in \sW^o_x$, we have $$(f^1_a \sL_1)_x =  \sR_x^* \otimes \Sym^3 a^{\perp}$$ under the identification $\sL_1 \cong \sR^* \otimes \Sym^3 \sW^*,$ where $$a^{\perp} := \{ h \in \sW_x^*,  \  h (a) =0\}.$$  Thus for any linearly independent $a_1, a_2 \in \sW^o_x$, we have
$$(f^1_{a_1} \sL_1)_x \cap (f^1_{a_2} \sL_1)_x =0,$$ which proves (iv).

For any $x \in M'$ and $a \in \sW^o_x$, we have $$(f^1_a \sV_{-1})_x =
\C a$$ under the identification $\sV_{-1} = \sW.$ Thus for any linearly independent $a_1, a_2 \in \sW^o_s$, we have $$(f^1_{a_1} \sV_{-1})_x \cap (f^1_{a_2} \sV_{-1})_x = 0,$$ which proves (vi).

For any $x \in M'$ and $a \in \sW^o_x$, we have $$(f^1_a \sV_{0})_x =
\sR_x^* \otimes \Sym^2 a^{\perp}$$ under the identification $\sV_{o} = \sR^* \otimes \Sym^2 \sW^*,$ where $$a^{\perp} := \{ h \in \sW_x^*,  \  h (a) =0\}.$$ Thus for any linearly independent $a_1, a_2 \in \sW^o_s$, we have $$(f^1_{a_1} \sV_{0})_x \cap (f^1_{a_2} \sV_{0})_x = 0,$$ which proves (vii).
\end{proof}

\begin{proposition}\label{p.generation}
In Assumption \ref{assumption},
let $\sP$ be the vector bundle defined as the kernel of the homomorphism $\wedge^2 \sL_{-1} \to \sL_{-2}$. We have the following properties of vector bundles on $M$.
\begin{itemize}
\item[(i)] $\sL_{-2}$ is $m$-tuply $\sO(1)$-generated for all $m\geq 1$.
\item[(ii)] $\sL_{-1}$ is $\sO(2)$-generated.
\item[(iii)] $\sV_{-1}$ is $\sO(1)$-generated.
\item[(iv)] $\sP$ is $\sO(3)$-generated.
\item[(v)] $\wedge^2 \sV_{-1}$ is $m$-tuply $\sO(1)$-generated for all $m \geq 1$.
\item[(vi)] $\sL_{-1} \otimes \sL_{-2}$ is $\sO(3)$-generated.
\item[(vii)] $\sV_{-1} \otimes \sL_{-2}$ is $\sO(2)$-generated and $m$-tuply $\sO(1)$-generated for all $m$.
\item[(viii)] $\sV_{-1} \otimes \sL_{-1}$ is $\sO(2)$-generated.
    \end{itemize}
\end{proposition}

\begin{proof}
Since $\sL_{-2}$ and $\wedge^2 \sV_{-1}$ are line bundles, (i) and (v) are direct consequence of Lemma \ref{c.split} and Lemma \ref{l.line}.
Also, (iii) is immediate from $\sV_{-1} = \sW$.

For any $x \in M'$ and any $a \in \sW^o_x$, we know from Assumption \ref{assumption} (3) and the identification $\sL_{-1} = \sR \otimes \Sym^3 \sW$ that $$(f^2_a \sL_{-1})_x = \sR_x \otimes a^3.$$ Thus $\sL_{-1}$ is $\sO(2)$-generated, proving (ii).

For any $x \in M'$ and any $a \in \sW^o_x$, Lemma \ref{c.split} shows $$f^3_a(\wedge^2 \sL_{-1}) \subset f_a^* \sP \subset f_a^* (\wedge^2 \sL_{-1}).$$
Under the identification $\sL_{-1} = \sR \otimes \Sym^3 \sW$, we have
$$(f^3_a \sP)_x = f^3_a(\wedge^2 \sL_{-1})_x = (\sR_x \otimes a^3) \wedge (\sR_x \otimes (a^2 \odot \sW_x).$$ Thus $\sP$ is $\sO(3)$-generated by Lemma \ref{l.Pspan}, proving (iv).

As $\sL_{-2}$ is a line bundle, (vi) is immediate from (ii) and Lemma \ref{c.split}.

As $\sL_{-2}$ is a line bundle, (iii) and Lemma \ref{c.split} implies that $\sV_{-1} \otimes \sL_{-2}$ is $\sO(2)$-generated. For any $x \in M'$ and $a \in \sW^o_,$
we have $$f_a^*(\sV_{-1} \otimes \sL_{-2}) \cong \sO(2) \oplus \sO(1)$$ from Lemma \ref{c.split}. Thus $$f^1_a(\sV_{-1} \otimes \sL_{-2})_x = (\sV_{-1} \otimes \sL_{-2})_x,$$ which shows $\sV_{-1} \otimes \sL_{-2}$ is $m$-tuply $\sO(1)$-generated for all $m \geq 1$. This proves (vii).

For any $x \in M'$ and $a \in \sW^o_x$, we have $$ f^2_a(\sV_{-1} \otimes \sL_{-1})_x \supset \sW_x \otimes \sR_x \otimes a^3$$ under the identification $\sV_{-1} \otimes \sL_{-1} = \sW \otimes \sR \otimes \Sym^3 \sW.$ Thus $\sV_{-1} \otimes \sL_{-1}$ is $\sO(2)$-generated, proving (viii). \end{proof}

We have the following vanishing result.

\begin{theorem}\label{t.vanishing}
In Assumption \ref{assumption}, for any $i, j <0$ and $k \geq i+j+1$, we have
$$\Hom( (\sV_i \oplus \sL_i) \wedge (\sV_j \oplus \sL_j), \sV_k \oplus \sL_k) = 0. $$\end{theorem}

The proof of Theorem \ref{t.vanishing} consists of the following series of lemmata
which check the vanishing of various components of $\Hom ((\sV_i \oplus \sL_i) \wedge (\sV_j \oplus \sL_j), \sV_k \oplus \sL_k).$

\begin{lemma}
For each integer $k$, we have $$\Hom(\wedge^2 \sL_{-2}, \sL_k)  = \Hom(\wedge^2 \sL_{-2}, \sV_k) =0.$$
\end{lemma}

\begin{proof}
This is automatic, since $\wedge^2\sL_{-2}$ is the zero bundle.
\end{proof}

\begin{lemma}
For each integer $k$, we have $$\Hom(\sL_{-1} \otimes \sL_{-2}, \sL_k) = \Hom (\sL_{-1} \otimes \sL_{-2}, \sV_k) =0.$$ \end{lemma}

\begin{proof} Since $\sL_{-1} \otimes \sL_{-2}$ is $\sO(3)$-generated, the vanishing follows from Lemma \ref{l.compare} and Proposition \ref{p.null}. \end{proof}

\begin{lemma} For any integers $i, j$ with $i \neq -2$, we have
$$\Hom(\wedge^2 \sL_{-1}, \sL_i) = \Hom (\wedge^2 \sL_{-1}, \sV_j) = 0.$$ \end{lemma}

\begin{proof}
From the short exact sequence $$ 0 \to \sP \to \wedge^2 \sL_{-1} \to \sL_{-2} \to 0,$$ it suffices to check the following for all integers $i,j$ with $i \neq -2$:
$$\Hom(\sP, \sL_k) = \Hom(\sP, \sV_k) = \Hom(\sL_{-2}, \sL_i) = \Hom (\sL_{-2}, V_j)  = 0.$$ This can be checked using the fact that $\sP$ is $\sO(3)$-generated and $\sL_{-2}$ is $4$-tuply $\sO(1)$-generated in Proposition \ref{p.generation}
by applying Lemma \ref{l.compare} and Proposition \ref{p.null}.
\end{proof}

\begin{lemma}
For any integers $i, j$ with $i \neq -2$, we have
$$\Hom(\wedge^2 \sV_{-1}, \sL_i) = \Hom (\wedge^2 \sV_{-1}, \sV_j) = 0.$$ \end{lemma}

\begin{proof}
Since $\wedge^2 \sV_{-1}$ is 4-tuply $\sO(1)$-generated by Proposition \ref{p.generation}, this follows from Lemma \ref{l.compare} and Proposition \ref{p.null}. \end{proof}

\begin{lemma}
For any integer $k$, we have
$$\Hom(\sV_{-1} \otimes \sL_{-2}, \sL_k) = \Hom(\sV_{-1} \otimes \sL_{-2}, \sV_k) = 0.$$
\end{lemma}

\begin{proof} Since $\sV_{-1} \otimes \sL_{-2}$ is $\sO(2)$-generated and 4-tuply $\sO(1)$-generated by Proposition \ref{p.generation}, this follows from Lemma \ref{l.compare} and Proposition \ref{p.null}. \end{proof}

\begin{lemma}
For any integer $k$, we have
$$\Hom(\sV_{-1} \otimes \sL_{-1}, \sL_k) = \Hom(\sV_{-1} \otimes \sL_{-1}, \sV_k) = 0.$$ \end{lemma}

\begin{proof}
Since $\sV_{-1} \otimes \sL_{-1}$ is $\sO(2)$-generated by Proposition \ref{p.generation}, the vanishing of all terms, except $\Hom(\sV_{-1} \otimes \sL_{-1}, \sL_{-1})$,   follows easily from Lemma \ref{l.compare} and Proposition \ref{p.null}.

For any $\sigma \in \Hom(\sV_{-1} \otimes \sL_{-1}, \sL_{-1}), x \in M'$ and $a \in \sW^o_x$, we have \begin{eqnarray*} f_a^* \sV_{-1} &\cong&  \sO(1) \oplus \sO \mbox{ and }  \\
f_a^* \sL_{-1} & \cong & \sO(2) \oplus \sO(1) \oplus \sO \oplus \sO(-1). \end{eqnarray*}
It follows that $f^*_a \sigma$ sends $(f_a^* \sV_{-1}) \otimes (f_a^2 \sL_{-1}) $ into $f^2_a \sL_{-1}$. The fiber of $f_a^2 \sL_{-1}$ at $x$ is $\sR_x \otimes a^3$ under the identification $\sL_{-1} = \sR \otimes \Sym^3 \sW$. Thus we have $$\sigma_x( (\sV_{-1})_x \otimes \sR_x \otimes a^3) \subset \sR_x \otimes a^3.$$ Since this holds for all $a \in \sW^o_x$, we have
$$\sigma_x( (\sV_{-1})_x \otimes  \xi) \subset  \xi$$
for all $\xi \in (\sL_{-1})_x.$
It follows that under the decomposition
$$ \Hom(\sV_{-1} \otimes \sL_{-1}, \sL_{-1}) = \Hom(\sV{-1}, {\rm End}^0(\sL_{-1})) \oplus \Hom(\sV_{-1}, \sO_M),$$ the element $\sigma$ belongs to $\Hom(\sV_{-1}, \sO_M)$. But since $\sV_{-1}$ is $\sO(1)$-generated by Proposition \ref{p.generation}, we have $\Hom(\sV_{-1}, \sO_M) =0$ from Lemma \ref{l.compare}. Thus $\sigma =0$. \end{proof}

This completes the proof of Theorem \ref{t.vanishing}.
We have the following corollary.

\begin{corollary}\label{c.vanishing}
In the setting of Theorem \ref{t.vanishing},
we have $$H^0(M, C^{\ell, 2}(\sA)/\partial (C^{\ell, 1}(\sA))) =0$$ for all $\ell \geq 1.$\end{corollary}

\begin{proof}
Let $U$ be a $\bG_0$-module. Recall (Definition 5.2 of \cite{HL}) that a filtration $0 = U^1 \subset U^2 \subset \cdots \subset U^k = U$ is $\bG_0$-reductive if the unipotent radical of $\bG_0$ sends $U^{i+1}$ to $U^i$ for each $i$.  Then the graded object $U^{i+1}/U^i$ of the filtration can be considered as a module of the reductive quotient of $\bG_0$ by its unipotent radical. In this case, the vanishing of $H^0(M, \sA \times^{\bG_0} U^{i+1}/U^i)$ for all $i$ implies the vanishing of $H^0(M, \sA \times^{\bG_o} U)$ (by Lemma 5.5 of \cite{HL}).

The $\bG_0$-reductive filtration $V_k \subset \fg_k$ from Lemma \ref{l.filt} gives rise to  $\bG_0$-reductive filtrations of $$C^{\ell,2}(\fg_0)/\partial(C^{\ell,1}(\fg_0)) \mbox{ and } C^{\ell, 2}(\fg_0) = \oplus_{i,j \in \N} \Hom( \fg_{-i} \wedge \fg_{-j}, \fg_{-i-j+ \ell}).$$ All graded objects of these filtrations appear as direct summands of the vector bundles in Theorem \ref{t.vanishing}. Thus Theorem \ref{t.vanishing} implies the vanishing of $H^0(M, C^{\ell, 2}(\sA)/\partial (C^{\ell, 1}(\sA))).$  \end{proof}

\section{Varieties of minimal rational tangents and the behavior of their fundamental forms}\label{s.vmrt}

Let us recall the notion of the variety of minimal rational
tangents. See \cite{HM99} and \cite{Hw01} for introductory
surveys.

\begin{definition}\label{d.VMRT}
Let $X$ be a uniruled projective manifold. For an irreducible
component  $\sK$ of the space of rational curves on $X$, denote by
$\rho: {\rm Univ}_{\sK} \to \sK$ and $\mu: {\rm Univ}_{\sK} \to X$ the universal family
morphisms. The component $\sK$ is called a {\em minimal dominating
component}, if for a general $x \in X$, the fiber $\mu^{-1}(x)$ is
non-empty and complete. Replacing ${\rm Univ}_{\sK}$ by its normalization if
necessary,  we can assume that $\rho$ is a $\BP^1$-bundle and
$\sK_x := \mu^{-1}(x)$ is a smooth projective variety for a general
$x \in X$ (e.g. Theorem 1.3 of \cite{Hw01}).    The {\em tangent map} $\tau: {\rm Univ}_{\sK} \dasharrow \BP TX$
associating  a smooth point of a rational curve to its tangent
direction is a rational map. By Theorem 3.4 of
\cite{Ke},  the tangent map
induces a morphism $\tau_x: \sK_x \to \BP T_x X$ for a general point $x \in X$, which is finite
over its image $\sC_x \subset \BP T_x X$. This image $\sC_x$ is the {\em variety of minimal
rational tangents} (VMRT) at $x$. Denote by $\sC \subset \BP
TX$ the proper image of $\tau$. Throughout, we will assume that $\sC \neq \BP TX$.
\end{definition}

\begin{proposition}\label{p.standard}
For a general point $\kappa \in \sK$,  denote the morphism
$$\mu|_{\rho^{-1}(\kappa)}: \rho^{-1}(\kappa) \cong \BP^1
\longrightarrow X$$ by $f: \BP^1 \to X$
 and the map $$\tau|_{\rho^{-1}(\kappa)}: \rho^{-1}(\kappa) \longrightarrow \BP TX$$ by $h: \BP^1 \to \BP TX$.
 Then the following holds.  \begin{enumerate} \item[(i)]
$f$ is an immersion and $$f^*TX \cong \sO(2) \oplus \sO(1)^p
\oplus \sO^{n-1-p}$$ where $\sO(2) \cong df(T \BP^1) \subset
f^*TX$, $p = \dim \sC_x$ and $n = \dim X.$ \item[(ii)] There exists a
neighborhood $U \subset {\rm Univ}_{\sK}$ of $\rho^{-1}(\kappa)$ in Euclidean topology such that  $\tau$ is a morphism on $U$, its image $\sY : = \tau(U) \subset \sC$
is a complex manifold and the morphism $\varpi: \sY \to \pi(\sY)$
induced by the natural projection $\pi: \BP TX \to X$ is a
submersion with connected fibers. \item[(iii)] Denoting by $T^{\pi}:= {\rm Ker}(d \pi)$  the relative
tangent bundle of $\pi: \BP TX \to X$, we have \begin{eqnarray*} h^* T^{\pi} &=& \Hom \left(T\BP^1,  f^*
TX/(df(T\BP^1))\right) \\ &\cong& \Hom(\sO(2), \sO(1)^p \oplus
\sO^{n-1-p}) \end{eqnarray*} and denoting by $T^{\varpi}$  the relative tangent
bundle of $\varpi: \sY \to \pi(\sY)$, we have
$$h^* T^{\varpi} = \Hom(\sO(2), \sO(1)^p) \cong \sO(-1)^p.$$
Consequently, the normal bundle $N$ of $\sY \subset \BP TX$
satisfies $$h^* N \cong \sO(-2)^{n-1-p}.$$ \item[(iv)] Given a
subvariety $X' \subset X$ of codimension $\geq 2$, we can choose a general
$\kappa \in \sK$ such that $f (\BP^1) \cap X' = \emptyset.$
\end{enumerate} \end{proposition}

\begin{proof} All are well-known.  (i) is from Theorem 1.2 of \cite{Hw01} and
(ii) is from Proposition 1.4 of \cite{Hw01}. (iii) is immediate from  Proposition 2.3 of \cite{Hw01}
and  the Euler sequence for the tangent bundle of the projective space.  (iv) is Proposition II.3.7 of \cite{Ko} (also Lemma 2.1 of \cite{Hw01}). \end{proof}

The following is also standard. See e.g. Proposition 1.2.1 of
\cite{HM99} or Proposition 2.4 of \cite{Hw01}.

\begin{proposition}\label{p.detP}
In Definition \ref{d.VMRT}, suppose that there exists
a vector subbundle $D \subset TX^o$ (i.e. a distribution $D$) on a Zariski open subset  $X^o \subset X$ such that $\sC_x
\subset \BP D_x$ at a general point $x \in X$. Let $\varphi_x:
\wedge^2 D_x \to T_x X/D_x$ be the
homomorphism defined by Lie brackets of local vector fields belonging to $D$. Then $\varphi_x(\wedge^2 P) =0$ for any 2-dimensional subspace $P \subset D_x$  tangent
to the affine cone $\widehat{\sC}_x \subset D_x$ of $\sC_x \subset \BP D_x$. \end{proposition}

Proposition \ref{p.detP} is used frequently in combination with the following result
when $X$ has Picard number 1. See  Proposition 1.2.2 of \cite{HM99} or Proposition 2.2 of
\cite{Hw01} for a proof.

\begin{proposition}\label{p.nonintegrable}
In Proposition \ref{p.detP}, if $X$ has Picard number 1, then $\varphi_x \neq 0$ for a general $x \in X$, i.e., the distribution $D$ is not integrable.
\end{proposition}

The following two propositions are from Proposition 2.2 and
Proposition 3.1 of \cite{Mk}. Although \cite{Mk} formulated them
only for certain special cases, the proof there can be  extended  easily to the
general setting stated below. We  provide the proofs here for the
readers' convenience.

\begin{proposition}\label{p.II}
Let $f: \BP^1 \to X$,  $h: \BP^1 \to \sY \subset \sC$ and $\varpi: \sY \to \pi(\sY) \subset X$ be as in
Proposition \ref{p.standard}. For $t \in \BP_1$, denote by  $Y_t$
the complex submanifold $\varpi^{-1}(f(t)) \subset \sC_{f(t)} \subset \BP
T_{f(t)}X$. Then the following holds. \begin{enumerate}
\item[(i)] For any two points $t, t' \in \BP^1$, the second
fundamental forms $\II_{h(t)}$ of $Y_t \subset \BP T_{f(t)}X$ at
$h(t)$ and $\II_{h(t')}$  of $Y_{t'} \subset \BP T_{f(t')}X$ at
$h(t')$ are isomorphic as systems of quadrics. In particular, by
the generality of the choice of $\kappa \in \sK$, we may choose  $\sY$ in Proposition \ref{p.standard}  such that the image $
{\rm Im}(\II_y)$ has the same dimension for any   $y \in \sY$,
defining a subbundle $N^{(2)}$ of the normal bundle $N$ of $\sY$ in $\BP TX$. \item[(ii)] When the
rank of $N^{(2)}$ is $r$, we have $h^* N^{(2)} \cong \sO(-2)^r$.
\item[(iii)] Let $\gamma: h^* T^{\pi} \to h^* N$ be the natural
projection induced by $h^* T(\BP TX) \to h^*N$.  By the
isomorphism $$h^* T^{\pi} = \Hom \left(T\BP^1, f^*
TX/(df(T\BP^1)) \right)$$ of Proposition \ref{p.standard} (iii),
the second osculating spaces of $\sY \subset \BP T X$ define a subbundle $\sT^2   \subset f^* TX$ with $\sT^1:=\ker\gamma\subset\sT^2$ such that
\begin{eqnarray*}
&& h^*T^\varpi=\Hom(T\BP^1, \sT^1/df(T\BP^1)),  \\
&& \gamma^{-1} (h^*N^{(2)}) = \rH( T\BP^1, \sT^2/\sT^1).
\end{eqnarray*}
Then
\begin{eqnarray*}
&& \sT^1=\sO(2)\oplus\sO(1)^p,  \\
&& \sT^2 \cong \sO(2) \oplus \sO(1)^p \oplus\sO^r, \\
&& f^*TX/\sT^2 \cong \sO^{n-1-p-r}.
\end{eqnarray*}
\end{enumerate}
\end{proposition}

\begin{proof}
Note that $\II_{h(t)}$ defines a holomorphic
section of  the vector bundle $$\Hom({\rm Sym}^2 h^*T^{\varpi}, h^*N)$$ on $\BP^1$. By Proposition \ref{p.standard} (iii), it
is a section of \begin{equation}\label{e.sym} \Hom({\rm Sym}^2 \sO(-1)^p, \sO(-2)^{n-1-p}) \cong \Hom({\rm Sym}^2 \sO^p, \sO^{n-1-p})\end{equation} which is a trivial vector bundle on $\BP^1.$
Hence for any two points $t,t'\in\BP^1$, the systems of quadrics $\II_{h(t)}$ and $\II_{h(t')}$ are isomorphic to each other, proving (i). Moreover, the images of $\II_{h(t)}$ give a trivial vector subbundle of $\sO^{n-1-p}$ in the right hand side of (\ref{e.sym}), or a vector subbundle isomorphic to $\sO(-2)^r$ in $\sO(-2)^{n-1-p}$ in the left hand side of (\ref{e.sym}), proving (ii).
(iii) follows from (ii) and Proposition \ref{p.standard} (iii).
\end{proof}

\begin{proposition}\label{p.III} In the setting of Proposition \ref{p.II}, assume
that  $\sC_x$ is irreducible at a general point $x \in X$ and its
affine cone $\widehat{\sC}_x$ spans a subspace $D_x \subset
T_x X$ of codimension 1, defining a distribution $D \subset TX$ outside a subset
${\rm Sing}(D) \subset X$ of codimension $\geq 2$. Assume furthermore that
\begin{enumerate} \item[(1)] the distribution $D$ is not
integrable, i.e., the homomorphism $\varphi_x $ of Proposition \ref{p.detP} is surjective for a general $x \in X$; \item[(2)] $n-p-r =3$ in the notation of Proposition
\ref{p.II}(ii), i.e., the subspace $N_y^{(2)} \subset N_y$ has codimension 2 for a general point $y \in \sY$.
\end{enumerate} Choose a
general $\kappa \in \sK$ such that $f: \BP^1 \to X$ satisfies $f(\BP^1) \cap {\rm
Sing}(D) = \emptyset,$  which is possible from Proposition \ref{p.standard} (iv).   Then for any two $t, t' \in \BP^1$, the third
fundamental forms $\III_{h(t)}$ of $Y_t \subset \BP T_{f(t)}X$ at
$h(t)$ and $\III_{h(t')}$  of $Y_{t'} \subset \BP T_{f(t')}X$ at
$h(t')$ are isomorphic.
Moreover, $$f^*(TX/D) \cong \sO(1) \mbox{ and
} f^*D \cong \sO(2) \oplus \sO(1)^p \oplus \sO^{n-p-3} \oplus
\sO(-1).$$
\end{proposition}

\begin{proof} Let $\sT^2 \subset f^*TX$ be as in Proposition
\ref{p.II} (iii). Then $\sT^2
\subset f^* D$ because $\sC_x$ is contained in $\BP D_x$.  By the condition (2),  the quotient bundle $f^*TX/\sT^2$ has rank 2
and $(f^*D)/\sT^2$ is a line subbundle of $f^*TX/\sT^2 \cong \sO^2.$
Thus $(f^*D)/\sT^2 \cong \sO(-k)$ for some $k \geq 0$. Since $\sT^2 \cong
\sO(2) \oplus \sO(1)^p \oplus \sO^{r}$ by Proposition \ref{p.II}
(iii),  the exact sequence $$0 \to \sT^2 \to f^*D \to (f^*D)/\sT^2 \to 0$$
splits, which implies
$$f^*D \cong \sO(2) \oplus \sO(1)^p \oplus \sO^{n-p-3} \oplus
\sO(-k) \mbox{ and } f^*TX/f^*D \cong \sO(k).$$

If $k=0$, we have a contradiction with the condition (1). To see this, let $\varphi: \wedge^2 D \to TX/D$ be the homomorphism given by Lie brackets of local vector fields belonging to the distribution $D$ outside ${\rm Sing}(D)$.
Under the assumption that $k =0$, $$ f^*D \cong \sO(2) \oplus
\sO(1)^p \oplus \sO^{n-p-2} \mbox{ and } f^*TX/f^*D \cong \sO.$$
 If $\varphi$ is not
identically zero, then for a general $x \in X$ and a general $\alpha \in \sC_x$ with the corresponding
1-dimensional subspace $\hat{\alpha} \subset
\hat{\sC}_x,$  we have
\begin{equation}\label{e.nonzero} \varphi(\hat{\alpha} \wedge  D_x) \neq 0 \end{equation}
because $\sC_x$ spans $D_x$.
 But we can
choose a general $\kappa \in \sK$ such that $h(t) = \alpha$. The
restriction of $f^*\varphi$ to $df(T\BP^1) \wedge f^*D \subset
f^*(\wedge^2 D)$ is a section of
$${\rm Hom}(T \BP^1 \otimes f^*D, \sO) = \Hom(\sO(2) \otimes
(\sO(2) \oplus \sO(1)^p \oplus \sO^{n-2-p}), \sO),$$ a
contradiction to $(\ref{e.nonzero}).$

Now assuming $k \geq 1$,  consider the third fundamental form
$$\III_{h(t)}: {\rm Sym}^3  T_{h(t)}Y_t \to h^*(N/N^{(2)})_t,$$
which factors through $$\III_{h(t)}: {\rm Sym}^3  T_{h(t)}Y_t \to
\Hom(T\BP^1,  (f^* D)/\sT^2)_{t}. $$ From Proposition \ref{p.standard}
(iii) and Proposition \ref{p.II} (iii), this is  a non-zero holomorphic
section of the bundle $\Hom( {\rm Sym}^3 (\sO(-1)^p), \sO(-k-2))$.
Thus $k=1$ and $\III_{h(t)}$ must be a constant section. Then
$$f^*(TX/D) \cong \sO(1) \mbox{ and
} f^*D \cong \sO(2) \oplus \sO(1)^p \oplus \sO^{n-p-3} \oplus
\sO(-1)$$ are immediate.
\end{proof}

\section{Proof of Theorems \ref{t.main}, \ref{t.deform} and \ref{t.rigid}}\label{s.proof}

To prove Theorem \ref{t.main}, we start with the following proposition relating the
condition of Theorem \ref{t.main} to Setup \ref{setup}.

\begin{proposition}\label{p.setup}
Let $X$ be a 7-dimensional uniruled projective manifold of Picard number 1 with a minimal dominating component $\sK$ such that
the VMRT $\sC_x \subset \BP T_x X \cong \BP^6$ at a general point $x \in X$ is projectively
isomorphic to $\BS \subset \BP \fg_{-1} \subset \BP \fg_{-}$ where $\fg_-= \fg_{-1} \oplus \fg_{-2}$ is as in Definition \ref{d.fg}.
Then there
exist \begin{itemize}
\item[(a)]  a Zariski open subset $M \subset X$; \item[(b)]  a vector subbundle $D\subset TM$ of rank $6$ inducing a filtration $F^{\bullet}= (F^{-1} = D \subset F^{-2} = TM)$
of type $\fg_-$; and \item[(c)] a $\bG_0$-structure $\sA \subset {\rm grFr}(F^{\bullet})$ subordinate to the filtration $F^{\bullet}$, \end{itemize} such that \begin{itemize} \item[(1)] a general member of $\sK$ is contained in $M$; \item[(2)] $M$ is simply connected; \item[(3)]  the VMRT $\sC_x $ is contained in $\BP D_x$ for each $x \in M$;  and \item[(4)] any element of $\sA_x \subset {\rm Isom}(\fg_{-}, {\rm Symb}_x(F^{\bullet}))$  sends $\BS \subset \BP \fg_{-1}$ to  the VMRT $\sC_x \subset \BP D_x$ for any $x \in M$. \end{itemize}  \end{proposition}

\begin{proof} From the assumption, the subset
 $M^+ \subset X$ defined by
$$M^+ :=  \{ x \in  X, \; (\sC_x \subset \BP T_x X)
\mbox{ is isomorphic to } (\BS
\subset \BP \fg_{-}) \}$$ is nonempty and  Zariski-open in $X$.
Denote by $D_x \subset T_x X$ the linear span of $\widehat{\sC}_x \subset T_s X$ for each $x \in M^+$, which is a hyperplane in  $T_x X.$ This defines a
distribution $D \subset T M^+$.  By Proposition \ref{p.nonintegrable}, the associated homomorphism $\varphi_x: \wedge^2 D_x \to T_x X/D_x$ is nonzero for a general $x \in M^+.$ Thus $$M:= \{ x \in M^+, \varphi_x \neq 0\}$$ is a nonempty Zariski-open subset in $X$.
By abuse of notation, write $D$ for the restriction $D|_M$ and let  $F^{\bullet}= (F^{-1} = D \subset F^{-2} = TM^+)$ be the associated filtration on $M$.

To prove (1), we claim first that a general member of $\sK$ lies in $M^+$.
It is clear that for a general member $f: \BP^1 \to X$ of $\sK$, we have
$f(\BP^1) \cap M^+ \neq \emptyset$.
Since the subvariety $\sC \subset \BP TX$ in Definition \ref{d.VMRT} is an irreducible subvariety, there exists a subvariety $\sZ \subset X$ of codimension at least 2 such that  the fiber $\sC_x \subset \BP T_xX$ of $\sC \to X$ at each $x \in X \setminus \sZ$ is an algebraic cycle of dimension 2 and  degree 4.
By Proposition \ref{p.standard} (iv),  we can choose a general $f$ such that $f(\BP^1) \cap \sZ = \emptyset$. Thus we may  assume that
$f^* \sC \subset  f^* \BP TX$ forms a well-defined family of proper algebraic cycles of dimension 2 and degree 4 over $\BP^1$
in the sense of Definition I.3.11 of \cite{Ko}. Then each irreducible component of the fiber $\sC_{f(t)}$ for any $ t \in \BP^1$ is a surface of degree $\leq 4$ in
$\BP^5$. But there exists  a smooth point $y$ of $\sC_{f(t)}$ such that $\II_{y}$ and $\III_{y}$ are isomorphic
to those of $S(1,3) \cong \BS$ by Proposition \ref{p.II} (i) and Proposition
\ref{p.III}.  Thus the linear span of $\sC_{f(t)}$ is of dimension 5 and $\sC_{f(t)}$ is isomorphic to $S(1,3)$ or $S(2,2)$ by Corollary \ref{c.DelPezzo}.
Since  $\sC_x$ is isomorphic to $S(1,3)$ for $x \in f(\BP^1) \cap M^+,$ we conclude that $\sC_{f(t)}$ is isomorphic to
$S(1,3)$ for any $ t\in \BP^1$ by Corollary \ref{c.deform}.  This proves the claim.

Now to prove $f(\BP^1) \subset M$ for a general member $f$ of $\sK$, assume the contrary, i.e. assume that $\varphi_{f(t_1)} =0$ at some point $f(t_1) \in f(\BP^1)$ of a general member $f: \BP^1 \to M^+$ of $\sK$.
 We can choose a section $v$ of the subsheaf $df(T\BP^1) \subset f^* D$ which vanishes at two points $t_2 \neq t_3$ different from $t_1$  such that $\varphi(v, \cdot)$ is a
 non-zero section of $f^*(D^* \otimes (T M^+/D))$. But $f^*\varphi(v, \cdot)$ has three
 distinct zeros at $t_1, t_2, t_3$. This is a contradiction to the isomorphism $$f^*(D^*
 \otimes (TM^+/D)) \cong \sO(-1) \oplus \sO^p \oplus \sO(1)^{n-p-3} \oplus \sO(2) $$ with $n=7$ and $ p=2, $ from Proposition \ref{p.III}. This proves (1).

(2) is immediate from (1).
Since $X$ is Fano, it is simply connected. As we have a complete curve on $M$
and $X$ has Picard number 1, the complement $X\setminus M$ has codimension at least 2 in $X$. It follows that $M$ is simply connected.

 By Proposition \ref{p.detP}, the Lie bracket homomorphism $\varphi: \wedge^2 D \to TM/D$ annihilates $\wedge^2 P$ for each  plane $P \subset D_x$ tangent to the affine cone $\widehat{\sC}_x \subset D_x$ of $\sC_x$ at every $x \in M$.
Then Proposition \ref{p.span} implies that any linear isomorphism $\fg_{-1} \to D_x$ sending $\BS$ to $\sC_x$ can be extended uniquely to a graded Lie algebra isomorphism $\fg_- \to {\rm Symb}_x (F^{\bullet})$ for a general $x \in M$.
Thus the set of isomorphisms from $\fg_{-1}$ to $D_x$ sending $\BS$ to $\sC_x$ determines a fiber subbundle $\sA \subset {\rm grFr}(F^{\bullet})$.  Since $M$ is simply connected, we may replace $\sA$ by a connected component of $\sA$ to assume that the fibers of $\sA \to M$ is connected. By   Proposition \ref{p.autS}, this is a $\bG_0$-structure subordinate to $F^{\bullet}$.
Properties (3) and (4) are automatic.
 \end{proof}

\begin{lemma}\label{l.split}
Since the situation of Proposition \ref{p.setup} belongs to Setup \ref{setup}, we have the associated vector bundles from Setup \ref{setup}. We have  a natural inclusion $\BP \sV_{-1,x} = \BP \sW  \subset \sC_x$ for each $x \in M$ coming from the natural inclusion $V_{-1} = W \subset \widehat{\BS}$.
Let $f: \BP^1 \to M$ be a general member of $\sK$ satisfying the properties of Proposition \ref{p.standard}. We have the vector bundles $\sT^i, i=0,1,2,3$ on $\BP^1$ determined by the osculating spaces of $\sY \subset \BP TM$ along $h:\BP^1 \to \sY.$
By the generality of $f$, we can assume that $$ h(o) = [{\rm d} f(T_o \BP^1)] \in \sC_x \setminus \BP \sV_{-1, x}.$$
Then we have the following.
\begin{itemize}
\item[(1)] The intersection $\sT^0 \cap f^* \sV_{-1}$ is zero.
\item[(2)] We have isomorphisms of vector bundles on $\BP^1$ $$\sT^0 \cong \sO(2), \sT^1 \cong \sO(2) \oplus \sO(1)^{\oplus 2}, \sT^2 \cong \sO(2) \oplus \sO(1)^{\oplus 2} \oplus \sO^{\oplus 2}.$$
    \item[(3)] We have an inclusion $f^* \sV_{-1}\subset\sT^2 $ as subbundles of  $f^* TM$ and the intersection $\sT^1 \cap f^*\sV_{-1}$ is a line subbundle, to be denoted by $\sW_{\flat}$.
            \item[(4)] Writing $\sW_{\sharp} := f^* \sW/\sW_{\flat}$, we have isomorphisms of line  bundles on $\BP^1$:
                $$\sT^0 \cong f^* \sR \otimes \Sym^3 \sW_{\flat}, \ \sN^3 := \sT^3 /\sT^2 \cong f^* \sR \otimes \Sym^3 \sW_{\sharp}. $$
\item[(5)] We have short exact sequences of vector bundles on $\BP^1$:
$$
0 \to \sW_{\flat} \to \sN^1 \to f^* \sR \otimes (\sW_{\sharp} \odot \Sym^2 \sW_{\flat}) \to 0, $$
$$0 \to \sW_{\sharp} \to \sN^2 \to f^* \sR \otimes (\sW_{\flat} \cdot \Sym^2 \sW_{\sharp}) \to 0.$$ \end{itemize} \end{lemma}

\begin{proof}
To prove (1), we need to show that for each $t \in \BP^1,$ the point $h(t) \in \sC_{f(t)}$ corresponds to a point on the open orbit $\BS \setminus \BP W$ of $\BS$ under $\bG_0$ in Proposition \ref{p.autS}. By assumption, this holds at $t = o \in \BP^1$. Then this holds for any other $ t \in \BP^1$ by Proposition \ref{p.II} (i), because the second fundamental form of $\BS$ at a point on the open orbit is not isomorphic to the one at a point of the closed orbit by Lemma \ref{l.IIz}.

(2) is from Proposition \ref{p.II} (iii).
(3) follows from (1) and Lemma \ref{l.ST}.
The short exact sequence of vector bundles on $M$ associated with (\ref{e.G0}) in Proposition \ref{p.ST}
$$0 \to \sW \to D \to \sR \otimes \Sym^3 \sW \to 0  $$
induces (4) and (5) as in Proposition \ref{p.ST}.
\end{proof}

\begin{proposition}\label{p.split}
In Lemma \ref{l.split}, we have the following isomorphisms of line bundles on $\BP^1$
\begin{equation}\label{e.split} f^* \sR \cong \sO(-1), \sW_{\flat} \cong \sO(1), \sW_{\sharp} \cong \sO, \sN^3 \cong \sO(-1),\end{equation} and $f^*(TM/\sD) \cong \sO(1) \cong \sQ$.
Consequently, the short exact sequences in Lemma \ref{l.split} (5) split and $$f^* \sW \cong \sO(1) \oplus \sO, $$ $$ f^* \sD \cong \sO(2) \oplus \sO(1)^{\oplus 2} \oplus \sO^{\oplus 2} \oplus \sO(-1).$$
\end{proposition}

\begin{proof}
By Lemma \ref{l.split} (4) and (5), we have \begin{eqnarray*}
\deg \sT^0 &=& \deg f^* \sR + 3 \deg \sW_{\flat} \\
\deg \sN^1 &=& \deg f^* \sR + 3 \deg \sW_{\flat} + \deg \sW_{\sharp} \\
\deg \sN^2 &=& \deg f^* \sR + \deg \sW_{\flat} + 3 \deg \sW_{\sharp} \\
\deg \sN^3 &=& \deg f^*\sR + 3 \deg \sW_{\sharp}. \end{eqnarray*}
Using Lemma \ref{l.split} (2), we obtain the degrees of the line bundles
$$\deg f^*R = -1, \deg \sW_{\flat} = 1, \deg \sW_{\sharp} = 0 \mbox{ and } \deg \sN^3 = -1,$$
which proves (\ref{e.split}).
From \begin{eqnarray*} \deg f^* D  &=& \deg \sT^0 + \sum_{1 \leq i \leq 3} \deg \sN^i \mbox{ and } \\ \deg f^* TM & =& \deg (\sO(2) \oplus \sO(1)^{\oplus 2} \oplus \sO^{\oplus 4}) =4, \end{eqnarray*} we obtain $\deg f^*(TM/D) = 1$ and $f^*(TM/D) \cong \sO(1).$
Finally,
the exact sequence $$0 \to \sT^2 \to \sT^3 \to  \sN^3 \to 0$$
should split from $H^1(\BP^1, \sT^2 \otimes (\sN^3)^*) =0,$ which implies
$$f^* D \cong \sT^3 \cong \sO(2) \oplus \sO(1)^{\oplus 2} \oplus \sO^{\oplus 2} \oplus \sO(-1).$$
 \end{proof}

\begin{proposition}\label{p.assumption}
In Proposition \ref{p.split}, let $\sK' \subset  \sK$ be the Zariski open subset whose members satisfy the condition in Lemma \ref{l.split} and let $M' \subset M$ be the Zariski open subset $$M' :=  \{ f(o) \in M, \ f \mbox{ belongs to } \sK'\}.$$
Then for each $x \in M',$
 there exists a Euclidean open subset $\sW_x^o \subset \sW_x$ and a member $f_a: \BP^1 \to M$ of $\sK'$  associated to  each $a \in \sW_x^o$ such that \begin{itemize} \item[(1)] $f_a(o) = x$; \item[(2)]  $f_a^* \sW$ contains a unique line subbundle $f^1_a \sW$ of degree 1 on $\BP^1$; \item[(3)]  $f_a^*\sW \cong f^1_a \sW \oplus \sO_{\BP^1}$; and \item[(4)] the fiber of $f^1_a\sW$ at $x$ is $\C a$. \end{itemize} Consequently, the Zariski open subsets $M' \subset M$ of $X$ satisfy Assumption \ref{assumption}.
\end{proposition}

\begin{proof}
For $f: \BP^1 \to M$ with $f(o) =x \in M'$ belonging to $\sK'$,
there is a line subbundle $f^1\sW \subset f^*\sW$ with $f^* \sW \cong f^1 \sW \oplus \sO_{\BP^1}$ from Proposition \ref{p.split}. We need to show that  as we vary $f$ among members of $\sK'$ through $x$, the locus of the fibers $(f^1 \sW)_x$ includes a Euclidean open subset in $\sW_x.$

The line subbundle  $f^* \sR \otimes (f^1 \sW)^{\otimes 3}$
corresponds to the unique $\sO(2)$-factor in $$f^* \sL_{-1} \cong f^*( \sR \otimes \Sym^3 \sW) \cong \sO(2) \oplus \sO(1) \oplus \sO \oplus \sO(-1).$$
This must be equal to  the image of ${\rm d}f (T \BP^1) \subset f^* D$ under the natural projection $D \to \sL_{-1}$ associated with $\fg_{-1} \to \fl_{-1}$.
Consider  deformations $f_t$ of $f$ with $f_t(o) =x$. Their tangent directions at $x$ give the germ of the surface $\sC_x \subset \BP \sD_x.$ It follows that the fibers of the $\sO(2)$-factors in $f_t^* \sL_{-1}$ move nontrivially in $\BP (\sR_x \otimes \Sym^3 \sW_x)$. Thus the locus of $(f^1_t \sW)_x$ as $f_t$ varies must contain  a Euclidean open subset $\sW^o_x$ in $\sW_x$. This proves the proposition. \end{proof}

We have the following lemma in the setting of Theorem \ref{t.deform}.

\begin{lemma}\label{l.deform}
Let $\pi: \sX \to \Delta$ be a smooth projective morphism from an
8-dimensional complex manifold $\sX$ to the unit disc $\Delta
\subset \C$. If the fiber $X_t:=\pi^{-1}(t)$  is biregular to ${\rm Lines}(\Q^5)$
for each $t\in \Delta \setminus \{0\}$, then the following holds.
\begin{itemize}
\item[(i)] The central fiber
$X_0:= \pi^{-1}(0)$ is a Fano manifold of Picard number 1 and has a unique minimal dominating component $\sK$ whose members have  degree $1$ with respect to the ample generator of the Picard group.

\item[(ii)] For a general point $x\in X_0$, the polarized variety $(\sK_x, \tau_x^*\mathcal{O}(1))$ is isomorphic to $(S(1, 3), L_{13})$ or $(S(2,2), L_{22})$.

\item[(iii)] The VMRT $\sC_x \subset \BP T_x X_0 \cong \BP^6$ at a general point $ x \in X_0$  is isomorphic to either $S(1,3)$ or $S(2,2)$ in a hyperplane $\BP^5 \subset \BP^6$.
\end{itemize}
  \end{lemma}

\begin{proof}
Since $\pi$ is smooth and projective, we have $H^k(X_0, \Z)=H^k(X_t, \Z)$ for any $t\neq 0$. It follows that $\Pic(X_t)\cong H^2(X_t, \Z)=\Z$ for all $t \in \Delta.$  Let $\sA$ be the line bundle on $\sX$ whose restriction $\sA_t$ to $X_t$ is the ample generator of $\Pic(X_t)$ for each $t \in \Delta$.  Since the anticanonical bundle of $X_t$ is isomorphic to $\sA_t^{\otimes 4}$  for $t\neq 0$, the same holds for $t=0$.

The variety  $X_t \cong {\rm Lines}(\Q^5), t \neq 0,$ has a  minimal dominating component $\sK^t$ corresponding to the family of lines covering ${\rm Lines}(\Q^5)$ under the Pl\"ucker embedding (see Section 2 of \cite{Hw97}). Limits of these rational curves give rational curves on $X_0$ having $\sA_0$-degree 1 such that their loci cover $X_0$.  So they give minimal dominating components of $X_0$ whose members have $\sA_0$-degree 1.  Conversely, given a minimal dominating component on $X_0$ with members of $\sA_0$-degree 1, its members
can be deformed to $X_t, t \neq 0$ by Kodaira's stability \cite{Kd}, which must correspond to lines on ${\rm Lines}(\Q^5)$.

Let us show that there is a unique such minimal dominating component on $X_0$. Let $m$ be the number of such minimal dominating components.
For a general point $x \in X_0$, the subscheme $\sK^0_x$ of rational curves through $x$ having  $\sA_0$-degree 1  is  smooth projective (e.g. Theorem 1.3 of \cite{Hw01}) and has at least $m$ components.   Choose a section $\sigma: \Delta \to \sX$
with $\sigma(0) =x.$ Then $\sK^0_x$ is the limit of the family of smooth projective varieties $$\{\sK^t_{\sigma(t)} \subset \sK^t, t \in \Delta, t \neq 0\}$$  consisting of members of $\sK^t$ passing through $\sigma(t)$. It follows that  $\sK^0_{\sigma(0)}$ is  irreducible, because so is $\sK^t_{\sigma(t)}$  for any $t \neq 0.$ Thus $m=1$, verifying (i). (ii) can obtained by applying Proposition \ref{p.deform}(i) to the family $\sK^t_{\sigma(t)}$, $t\in\Delta$.

To prove (iii), it suffices to show that the image of $\tau_x: \sK_x \to \sC_x \subset \BP T_x X_0$ spans a subspace of dimension 6 in $T_x X_0$. Suppose not. Then the linear span of $\sC_x$ gives a distribution $D$ of rank $\leq 5$, which is not integrable by Proposition \ref{p.nonintegrable}. By Corollary 1 of \cite{HM04},  the image $\sC_x$ of $\tau_x$ is a surface obtained from the projection of $S(1,3)$ or $S(2,2)$ under a nontrivial projection of the form given in  Proposition \ref{p.span2} (ii) or (iv). Thus Proposition \ref{p.span2} and Proposition \ref{p.detP} imply that $D$ must be integrable, a contradiction. \end{proof}

Applying Lemma \ref{l.deform} to  Theorem \ref{t.PP}, we obtain the following, which can be also checked by a more concrete computation.

\begin{corollary}\label{c.bX}
There is a unique minimal dominating component ${\bf K}$  on $\bX \subset \BP(V(\omega_1) \oplus V(\omega_2))$ in Section \ref{s.I} whose members  have degree 1 with respect to the ample generator of ${\rm Pic}(\bX)$. The corresponding VMRT at a general point of $\bX$ is isomorphic to $S(1,3) \subset \BP^5 \subset \BP^6$. \end{corollary}

\begin{proof} By Theorem \ref{t.PP}, we can use Lemma \ref{l.deform} to see the first sentence.  Moreover, its VMRT at a general point is isomorphic to either $S(1,3)$ or $S(2,2)$.
But if it is isomorphic to $S(2,2)$,
Theorem \ref{t.Mok} implies that $\bX$ is isomorphic to ${\rm Lines}(\Q^5)$, a contradiction.  \end{proof}

We recall the following version of Cartan-Fubini type extension theorem, the proof of which is contained in \cite{HM01} (see Theorem 3.9 of \cite{Hw19} and its proof).

\begin{theorem}\label{t.CFO}
Let $X$ and $\widetilde{X}$ be Fano manifolds of Picard number 1.
 Let $\sK$ (resp. $\widetilde{\sK}$) be a minimal dominating component of the space of rational curves on $X$ (resp. $\widetilde{X}$). Assume \begin{itemize}
  \item[(i)] the VMRT $\sC_x$ and $\widetilde{\sC}_{\widetilde{x}}$ are irreducible for  general $x \in X$ and $\widetilde{x} \in \widetilde{X}$; and
      \item[(ii)]  general members of $\sK$ and $\widetilde{\sK}$ are nonsingular. \end{itemize}
          If there are members $A \subset X$ of $\sK$ and $\widetilde{A} \subset \widetilde{X}$ of $\widetilde{\sK}$ with open neighborhoods $A \subset U \subset X$ and $\widetilde{A} \subset \widetilde{U} \subset \widetilde{X}$ and  a biholomorphic map  $$\Phi: U \cong \widetilde{U} \mbox{ with } \Phi(A) = \widetilde{A},$$ then there exists a biregular morphism $ X \to \widetilde{X}$ which induces $\Phi$. \end{theorem}

Now we are ready to finish the proofs of Theorems \ref{t.main}, \ref{t.deform} and \ref{t.rigid} .

\begin{proof}[Proof of Theorem \ref{t.main}]
We can apply Proposition \ref{p.setup} and Proposition \ref{p.assumption} to obtain Zariski-open subset $M' \subset M \subset X$ satisfying Assumption \ref{assumption}.
Theorem \ref{t.HL} can be applied to $\sA$ on the simply connected manifold $M$ because of
Proposition \ref{p.kim} and Corollary \ref{c.vanishing}.
Thus we have an open immersion $\phi: M \to \bG/\bG^0$. The latter is a Zariski open subset in $\bX,$ the open orbit of the automorphism group of $\bX$.
 The images of members of $\sK$ are sent to rational curves in $\bX$ which have
 anticanonical degree $4$, i.e. have degree 1 with respect to the ample generator of ${\rm Pic}(\bX)$. It follows that $\phi$ sends general members of $\sK$ to
 general members of  ${\bf K}$  in Corollary \ref{c.bX}.   Thus $\phi$ can be extended to a biregular morphism $X \cong \bX$ by Theorem \ref{t.CFO}.
\end{proof}

\begin{proof}[Proof of Theorem \ref{t.deform}]
From Lemma \ref{l.deform}, the VMRT   $\sC_x $ at a general point $x\in X_0$ is isomorphic to either $S(2,2)$ or $S(1,3).$
 In the former case,  we have $X_0 \cong {\rm Lines}(\Q^5)$ by Theorem \ref{t.Mok}. In the latter case, we have $X_0 \cong \bX$ by Theorem \ref{t.main}.
 \end{proof}

\begin{proof}[Proof of Theorem \ref{t.rigid}]
Arguing in the same way as in the proof of Lemma \ref{l.deform}, applying
Proposition \ref{p.deform} (ii) instead of Proposition \ref{p.deform} (i), we see
 that the variety of minimal rational tangents $\sC_x$ at a general point $x \in X_0$ is
isomorphic to $S(1,3) \cong \BS \subset \BP^5 \subset \BP^6$.
 Thus  $X_0$
must be biholomorphic to $\bX$ by Theorem \ref{t.main}. \end{proof}

\bigskip
{\bf Acknowledgment}
We are grateful to the referees for valuable suggestions and pointing out an error in our use of Theorem \ref{t.DelPezzo} in the first version of the paper.

\bigskip

{\bf\Large Declaration}

\bigskip

{\bf Conflict of interest} The authors declare that there is no conflict of interest.

\bigskip
Jun-Muk Hwang(jmhwang@ibs.re.kr)

\smallskip

Center for Complex Geometry,
Institute for Basic Science (IBS),
Daejeon 34126, Republic of Korea

\bigskip

Qifeng Li(qifengli@sdu.edu.cn)

\smallskip

School of Mathematics,
Shandong University,
Jinan 250100, China


\begin{thebibliography}{KSWZ}
\bibitem[Br]{Br} Brieskorn, E.: \"Uber holomorphe $\BP_n$-B\"undel \"uber $\BP_1$. Math. Ann. {\bf 157} (1965) 343-357

\bibitem[EH]{EH} Eisenbud, D. and Harris, J.:  On varieties of minimal degree (a centennial account). {\em Algebraic geometry, Bowdoin, 1985 (Brunswick, Maine, 1985)}, 3--13, Proc. Sympos. Pure Math., 46, Part 1, Amer. Math. Soc., Providence, RI, 1987
\bibitem[FH]{FH} Fulton, W. and Harris, J.: {\em Representation theory. A first course}.
Graduate Texts in Mathematics. 129, Springer-Verlag, New York, 1991

\bibitem[Ha]{Ha} Hartshorne, R.: {\em Algebraic geometry}. Graduate Texts in Mathematics. 52, Springer-Verlag, New York, 1977

\bibitem[HH]{HH} Hong, J. and Hwang, J.-M.:
Characterization of the rational homogeneous space associated to a
long simple root by its variety of minimal rational tangents. {\em
Algebraic Geometry in East Asia -- Hanoi 2005}, Advanced Studies
in Pure Mathematics {\bf 50} (2008) 217-236
\bibitem[Hw97]{Hw97} Hwang, J.-M.: Rigidity of homogeneous contact manifolds under Fano
deformation. J. reine angew. Math. {\bf 486} (1997) 153-163
\bibitem[Hw01]{Hw01} Hwang, J.-M.: Geometry of minimal rational curves
on Fano manifolds. {\em School on Vanishing Theorems and Effective
Results in Algebraic Geometry (Trieste, 2000)}, 335--393, ICTP
Lect. Notes, 6, Abdus Salam Int. Cent. Theoret. Phys., Trieste,
2001

\bibitem[Hw19]{Hw19}  Hwang, J.-M.: An application of Cartan's equivalence method to Hirschowitz's conjecture on the formal principle. Ann. Math.   {\bf 189} (2019) 979--1000

\bibitem[HL]{HL} Hwang, J.-M. and Li, Q.: Characterizing symplectic Grassmannians by varieties of minimal rational tangents. J. Diff. Geom. {\bf 119} (2021) 309--381

\bibitem[HM99]{HM99} Hwang, J.-M. and Mok, N.: Varieties of minimal
rational tangents on uniruled projective manifolds. {\em Several
complex variables (Berkeley, CA, 1995--1996)}, 351--389, Math.
Sci. Res. Inst. Publ., 37, Cambridge Univ. Press, Cambridge, 1999
\bibitem[HM01]{HM01} Hwang, J.-M. and Mok, N.:
Cartan-Fubini type extension of holomorphic maps for Fano
manifolds of Picard number 1. Journal Math. Pures Appl. {\bf 80}
(2001) 563-575
\bibitem[HM02]{HM02} Hwang, J.-M. and Mok, N.: Deformation rigidity of the rational homogeneous
space associated to a long simple root.  Ann. scient. Ec. Norm.
Sup. {\bf 35} (2002) 173-184
\bibitem[HM04]{HM04}
Hwang, J.-M. and Mok, N.: Birationality of the tangent map
for minimal rational curves. Asian J. Math. {\bf 8} (2004) 51--63

\bibitem[IL]{IL} Ivey, T. and Landsberg, J. M.: {\em Cartan for
beginners: differential geometry via moving frames and exterior
differential systems.} Graduate Studies in Mathematics, 61.
American Mathematical Society, Providence,  2003

\bibitem[Ke]{Ke} Kebekus, S.: Families of singular rational curves. J. Alg. Geom. {\bf 11} (2002) 245-256
\bibitem[Kim]{Kim} Kim, S.: Geometric structures modeled on smooth projective horospherical varieties of Picard number one. Transformatin Groups {\bf 22} (2017) 361-386

\bibitem[Kd]{Kd} Kodaira, K.: On stability of compact submanifolds of complex manifolds.  Amer. J. Math. {\bf 85} (1963) 79-94

\bibitem[Ko]{Ko} Koll\'ar, J.: {\em Rational curves on algebraic varieties.}
Ergebnisse der Mathematik und ihrer Grenzgebiete, 3 Folge, Band
32, Springer Verlag,  1996



\bibitem[Mk]{Mk} Mok, N.: Recognizing certain rational homogeneous
manifolds of Picard number 1 from their varieties of minimal
rational
 tangents.  Third International Congress of Chinese Mathematicians.
 Part 1, 2, 41-61, AMS/IP Stud. Adv. Math., 42, pt.1, 2,
  Amer. Math. Soc., Providence, RI, 2008
\bibitem[Pa]{Pa} Pasquier, B.:  On some smooth projective two-orbit varieties with Picard number 1. Math. Annalen {\bf 344} (2009) 963-987
\bibitem[PP]{PP} Pasquier, B. and
Perrin, N.:  Local rigidity of quasi-regular varieties.
 Math. Z. {\bf 265} (2010) 589-600

\bibitem[Ya]{Ya} Yamaguchi, K.: Differential systems associated with simple graded Lie algebras.
Adv. Study Pure Math. {\bf 22}  (1993) 413-494
\end{thebibliography}
\end{document}